\documentclass[a4paper,11pt,twoside,reqno]{amsart}

\usepackage[utf8]{inputenc}
\usepackage[plainpages=false,pdfpagelabels=true]{hyperref}
\usepackage{amssymb,amsthm,slashed,dsfont}
\usepackage[margin=1in]{geometry}
\usepackage[all]{xy}
\usepackage{color}

\newtheorem{Satz}{Theorem}[section]
\newtheorem{Prop}[Satz]{Proposition}
\newtheorem{Lem}[Satz]{Lemma}

\theoremstyle{definition}

\newtheorem{Bem}[Satz]{Remark}
\newtheorem{Bsp}[Satz]{Example}
\newcommand{\D}{\slashed{D}}
\newcommand{\p}{\slashed{\partial}}
\newcommand{\tr}{\operatorname{tr}}
\newcommand{\sff}{\mathrm{I\!I}}

\newcommand{\scal}{\mathrm{scal}}
\newcommand{\ric}{\mathrm{Ric}}
\newcommand{\trace}{\operatorname{tr}}

\newcommand{\dv}{\text{ }dV}

\allowdisplaybreaks[1]

\renewcommand{\epsilon}{\varepsilon}

\newcommand{\R}{\ensuremath{\mathbb{R}}}

\newcommand{\N}{\ensuremath{\mathbb{N}}}

\numberwithin{equation}{section}

\title[Global existence of Dirac-wave maps with curvature term]{Global existence of Dirac-wave maps with curvature term on expanding spacetimes}
\author{Volker Branding}
\date{\today}
\address{University of Vienna, Faculty of Mathematics\\
Oskar-Morgenstern-Platz 1, 1090 Vienna, Austria\\}
\email{volker.branding@univie.ac.at}
\author{Klaus Kr\"oncke}
\address{University of Hamburg\\
Department of Mathematics
Bundesstraße 55,
20146 Hamburg, Germany}
\email[]{klaus.kroencke@uni-hamburg.de}
\begin{document}

\begin{abstract}
We prove the global existence of Dirac-wave maps with curvature term with small initial data on globally hyperbolic manifolds of arbitrary dimension
which satisfy a suitable growth condition.
In addition, we also prove a global existence result for wave maps under similar assumptions.
\end{abstract} 

\maketitle

\section{Introduction and results}
Wave maps are among the fundamental variational problems in differential geometry.
They are defined as critical points of the Dirichlet energy for a map between two manifolds,
where one assumes that the domain manifold is Lorentzian and the target manifold is Riemannian.
More precisely, let  $(M,h)$ be a globally hyperbolic Lorentzian manifold, $(P,G)$ a Riemannian manifold and $\phi\colon M\to P$ a map.
Squaring the norm of its differential gives rise
to the Dirichlet energy, whose critical points are given by the \emph{wave map equation}, which is
\(\tau(\phi)=0\), where \(\tau(\phi):=-h^{\alpha\beta}\nabla_{\partial_\alpha}d\phi(\partial_\beta)\).

The wave map equation is a second-order semilinear hyperbolic system, that is the natural analogue of
the harmonic map equation for maps between Riemannian manifolds.
Wave maps are also well-studied in the physics literature, they appear as critical points
of the \emph{Polyakov action} in bosonic string theory for a string with Lorentzian worldsheet.

There are many articles that study wave maps in the case that the 
domain is Minkowski space. We cannot give an exhaustive list of the results
here, but want to mention the influential works of Klainerman \cite{MR837683},
Tataru \cite{MR1641721}, Klainerman and Selberg \cite{MR1901147}, Klainerman and Machedon \cite{MR1381973},
Tao \cite{MR1820329,MR1869874} and Shatah and Struwe \cite{MR1674843,MR1890048}.

There are less articles that consider the case of a domain being a non-flat globally hyperbolic manifold.
Here, we want to mention the works of Choquet-Bruhat \cite{MR1843031,MR1776183} for wave maps 
on Robertson-Walker spacetimes and several recent articles that consider wave maps on non-flat backgrounds \cite{MR3493440,MR3441208,MR2984144,MR3082240}.
To obtain an overview on the current status of research on the wave map equation
we refer to the recent book \cite{MR3585834}.

In modern quantum field theory one considers extensions of the wave map system,
one of them being the supersymmetric nonlinear \(\sigma\)-model \cite{MR626710}.
Recently, there has been a lot of interested in this model from a mathematical perspective.
In mathematical terms, the model leads to a geometric variational problem that couples
a map between manifolds to spinor fields.

Most of the mathematical research on this model so far is concerned with the case that both manifolds
are Riemannian leading to the notion of Dirac-harmonic maps
\cite{MR2262709} and Dirac-harmonic map with curvature term \cite{MR3333092, MR2370260},
which represent semilinear elliptic problems.
At present, many results regarding the geometric and analytic structure of Dirac-harmonic maps and
Dirac-harmonic maps with curvature are known \cite{MR3558358}, but no existence result 
for these kind of equations could be achieved.

In case that the domain manifold has a Lorentzian metric the critical points of the supersymmetric nonlinear \(\sigma\)-model 
lead to a system of the wave map equation coupled to spinor fields.
For this system two existence results are available
\cite{branding2017energy,MR2138082} for the domain being two-dimensional Minkowski space. 

In order to couple the wave map equation to spinor fields we have to recall
some concepts from spin geometry on globally hyperbolic manifolds.
We have to make the additional assumption that the manifold \(M\) is spin, 
which guarantees the existence of the spinor bundle \(SM\). 
Sections in the spinor bundle are called spinors.
Moreover, we fix a spin structure. 
The spinor bundle is a vector bundle on which we choose a metric connection compatible
with the hermitian scalar product denoted by \(\langle\cdot,\cdot\rangle_{S M}\).
On the spinor bundle we have the Clifford multiplication of spinors with tangent vectors,
which satisfies
\begin{align*}
\langle\psi,X\cdot\xi\rangle_{SM}=\langle X\cdot\psi,\xi\rangle_{SM}
\end{align*}
for all \(X\in TM,\psi,\xi\in\Gamma(SM)\).
In addition, the Clifford relations
\begin{align*}
X\cdot Y+Y\cdot X =-2h(X,Y)
\end{align*}
hold for all \(X,Y\in TM\), where \(h\) represents the metric on \(M\).

The natural operator acting on spinors is the \emph{Dirac operator},
which is defined as the composition of applying the covariant derivative first
and Clifford multiplication in the second step. More precisely, the Dirac 
operator acting on spinors is given by
\begin{align*}
\p:=h^{\alpha\beta}\partial_\alpha\cdot\nabla_{\partial_\beta}.
\end{align*}
The Dirac operator is a linear first order hyperbolic differential operator.
For more background on spin geometry on globally hyperbolic manifolds we refer to
\cite{MR2121740,MR701244}.
The Dirac operator itself is anti-self-adjoint with respect to the \(L^2\)-norm.
However, the combination \(i\p\) yields a self-adjoint operator, that is
\begin{align*}
\int_M\langle i\p\xi,\psi\rangle \dv_h=\int_M\langle\xi,i\p\psi\rangle \dv_h.
\end{align*}
The square of the Dirac operator satisfies the \emph{Schrödinger-Lichnerowicz}
formula 
\begin{align*}
\p^2=\Box+\frac{\scal^M}{4},
\end{align*}
where \(\scal^M\) denotes the scalar curvature of the manifold \(M\). 
Note that \(\p^2\) is a wave-type operator.

In quantum field theory one usually considers spinors that are twisted by some
additional vector bundle. Here, we consider the case that the spinor bundle
is twisted with the pullback of the tangent bundle from the target manifold.
More precisely, we are considering the bundle \(SM\otimes\phi^\ast TP\)
and sections in this bundle are called \emph{vector spinors}.
We obtain a connection on \(SM\otimes\phi^\ast TP\) by setting
\begin{align*}
\nabla^{SM\otimes\phi^\ast TP}=\nabla^{SM}\otimes \mathds{1}^{\phi^\ast TP}+\mathds{1}^{SM}\otimes\nabla^{\phi^\ast TP}.
\end{align*}
Note that Clifford multiplication on the twisted bundle is defined by acting only on the first factor.
This allows us to define the Dirac operator acting on vector spinors as follows
\begin{align*}
\D:=h^{\alpha\beta}\partial_\alpha\cdot\nabla^{SM\otimes\phi^\ast TN}_{\partial_\beta}.
\end{align*}
We assume that the connection on \(\phi^\ast TP\) is metric and thus the operator
\(i\D\) is also self-adjoint with respect to the \(L^2\)-norm.
After these considerations we are ready to present the energy functional for Dirac-wave maps with curvature term
\begin{align}
\label{action-functional-quer}
S(\phi,\tilde\psi,\tilde h)=\frac{1}{2}\int_{M}(|d\phi|^2+\langle\tilde\psi,i\tilde\D\tilde\psi\rangle-\frac{1}{6}\langle\tilde\psi,R^P(\tilde\psi,\tilde\psi)\tilde\psi\rangle)\dv_{\tilde h}.
\end{align}
In the last term the indices are contracted as follows
\begin{align*}
\langle\tilde\psi,R^P(\tilde\psi,\tilde\psi)\tilde\psi\rangle_{SM\otimes\phi^\ast TP}
=R_{IJKL}\langle\tilde\psi^I,\tilde\psi^K\rangle_{SM}\langle\tilde\psi^J,\tilde\psi^L\rangle_{SM},
\end{align*}
which ensures that the energy functional is real-valued.
We would like to point out that in the physics literature (see e.g. \cite{MR626710}) one usually considers Grassmann-valued spinors in the analysis
of \eqref{action-functional-quer}. However, we want to use methods from the geometric calculus of variation
and due to this reason we are employing standard spinors.

The critical points of \eqref{action-functional-quer} can easily be calculated as
\begin{align}
\label{phi-dwmap-quer}
\tau(\phi)=&\frac{1}{2}h^{\alpha\beta}R^P(\tilde\psi,i\partial_\alpha\cdot\tilde\psi)d\phi(\partial_\beta)-\frac{1}{12}\langle(\nabla R^P)^\sharp(\tilde\psi,\tilde\psi)\tilde\psi,\tilde\psi\rangle, \\
\label{psi-dwmap-quer}
i\tilde\D\tilde\psi=&\frac{1}{3}R^P(\tilde\psi,\tilde\psi)\tilde\psi,
\end{align}
where \(\tau(\phi):=-h^{\alpha\beta}\nabla_{\partial_\alpha}d\phi(\partial_\beta)\) represents the wave map operator.
The solutions of \eqref{phi-dwmap-quer}, \eqref{psi-dwmap-quer} are called \emph{Dirac-wave maps with curvature term}.

We are able to provide the first existence result for the Dirac-wave map with curvature term system which is as follows:

\begin{Satz}\label{result_diracwavemap}
Let $\tilde{g}_t$ be a smooth family of complete Riemannian metrics on $\Sigma^{n-1}$, $N\in C^{\infty}(\R\times\Sigma)$ with $0< A\leq N\leq B<\infty$ and
 $(M^n,\tilde{h})=(\R\times\Sigma,-N^2dt^2+\tilde{g}_t)$ be a globally hyperbolic Lorentzian spin manifold that satisfies the following condition: 
 There exists a monotonically increasing smooth function $s:\R\to\R_+$ with $\int_0^\infty s^{-1} dt <\infty$, such that the conformal metric 
\begin{align*}
h=(Ns)^{-2}\tilde{h}=-s^{-2}dt^2+g_t
\end{align*}
has bounded geometry, the metrics $g_t$ have a uniform Sobolev constant and 
\begin{align}
\label{assumptions-diracwavemap}
\|N\|_{C^k(g_t)}+\|\nabla_{\nu}N\|_{C^k(g_t)}+\|\sff\|_{C^k(g_t)}+s\|\sff\|_{L^\infty}\leq C<\infty
\end{align}
for all $k\in\N$. Here, $\nu$ is the future-directed unit normal of the hypersurfaces $\left\{t\right\}\times\Sigma$ and $\sff$ is their second fundamental form.

Then if in addition, the Riemannian manifold $(P,G)$ has bounded geometry, there exists for each \(r\in\N\) with \(r>\frac{n-1}{2}\) an \(\epsilon>0\)
such that if the initial data $(\phi_0,\phi_1,\psi_0)$ for the system \eqref{phi-dwmap-quer}, \eqref{psi-dwmap-quer} satisfies
\begin{align}
\|\phi_0\|_{H^{r+1}(\tilde{g}_0)}+\|\phi_1\|_{H^r(\tilde{g}_0)}+\|\psi_0\|_{H^{r}(\tilde{g}_0)}<\epsilon,
\end{align}
the unique solution of the system \eqref{phi-dwmap-quer}, \eqref{psi-dwmap-quer} with initial data
\(\phi|_{t=0}=\phi_0,\partial_t\phi|_{t=0}=\phi_1,\psi|_{t=0}=\psi_0\) exists
for all times \(t\in [0,\infty)\) and satisfies 
\begin{align*}
&\phi\in C^0([0,\infty),H^{r+1}(\Sigma,P))\cap  C^1([0,\infty),H^r(\Sigma,P)), \\
&\psi\in C^0([0,\infty),H^{r}(M,SM\otimes\phi^\ast TP))\cap C^1([0,\infty),H^{r-1}(M,SM\otimes\phi^\ast TP)).
\end{align*}
\end{Satz}
\begin{Bem}
Under the same assumptions, the proof of Theorem \ref{result_diracwavemap} also implies global existence of \emph{Dirac-wave maps},
which satisfy an equation slightly simpler than \eqref{phi-dwmap-quer}, \eqref{psi-dwmap-quer}. 
In this case, the second term on the right hand side of \eqref{phi-dwmap-quer} and the right-hand side of \eqref{psi-dwmap-quer}  both vanish.
\end{Bem}
Along the line of Theorem \ref{result_diracwavemap} we obtain the following result for the wave map equation
generalizing the results from \cite{MR1843031,MR1776183}:

\begin{Satz}\label{result_wavemap}
Let $\tilde{g}_t$ be a smooth family of complete Riemannian metrics on $\Sigma^{n-1}$, $N\in C^{\infty}(\R\times\Sigma)$ with $0< A\leq N\leq B<\infty$ and
 $(M^n,\tilde{h})=(\R\times\Sigma,-N^2dt^2+\tilde{g}_t)$ be a globally hyperbolic Lorentzian manifold that satisfies the following condition: 
There exists a monotonically increasing smooth function $s:\R\to\R_+$ with $\int_0^\infty s^{-1} dt <\infty$, such that for the conformal metric 
\begin{align*}
h=(Ns)^{-2}\tilde{h}=-s^{-2}dt^2+g_t
\end{align*}
the metrics $g_t$ admit a uniform Sobolev constant and 
\begin{align}
\|R^{g_t}\|_{C^k(g_t)}+\|N\|_{C^k(g_t)}+\|\nabla_{\nu}N\|_{C^k(g_t)}+\|\sff\|_{C^k(g_t)}\leq C<\infty
\end{align}
for all $k\in\N$. Here, $\nu$ is the future-directed unit normal of the hypersurfaces $\left\{t\right\}\times\Sigma$ and $\sff$ is their second fundamental form.

Then if in addition, the Riemannian manifold $(P,G)$ has bounded geometry, there exists for each \(r\in\N\) with \(r>\frac{n-1}{2}\) an \(\epsilon>0\)
such that if the initial data $(\phi_0,\phi_1)$ for the wave map equation satisfies
\begin{align}
\|\phi_0\|_{H^{r+1}(\tilde{g}_0)}+\|\phi_1\|_{H^r(\tilde{g}_0)}<\epsilon,
\end{align}
 the unique wave map with initial data
\(\phi|_{t=0}=\phi_0,\partial_t\phi|_{t=0}=\phi_1\) exists
for all times \(t\in [0,\infty)\) and satisfies 
\begin{align*}
&\phi\in C^0([0,\infty),H^{r+1}(\Sigma,P))\cap C^1([0,\infty),H^r(\Sigma,P)).
\end{align*}
\end{Satz}

\begin{Bem}
If we compare the assumptions of Theorem \ref{result_diracwavemap} and Theorem \ref{result_wavemap} we observe 
that one also has to control the $R_{0ij0}$-components of the curvature tensor in the result for Dirac-wave maps
with curvature term. This curvature contribution appears when computing the curvature of the spinor bundle. 
This explains why we have to assume bounded geometry of the Lorentzian manifold in Theorem \ref{result_diracwavemap}
but only bounded geometry of the Riemannian slices in Theorem \ref{result_wavemap}.

In addition, we also have to demand a decay of the second fundamental form \eqref{assumptions-diracwavemap} in Theorem \ref{result_diracwavemap}
which originates from the choice of a positive definite scalar product on the spinor bundle that is no longer metric.
\end{Bem}
\begin{Bem}
It seems that the class of Lorentzian manifolds that we are considering is the appropriate setting that guarantees a nice long-time behavior for solutions of various nonlinear wave equations with small initial data.
We therefore think that our approach can be also used to establish existence results for a large class of other second order hyperbolic systems arising in mathematical physics.
\end{Bem}
\begin{Bsp}
There are many spacetimes that satisfy the assumptions in Theorem \ref{result_diracwavemap} and Theorem \ref{result_wavemap}. 
The simplest class is the class of Robertson-Walker spacetimes $-dt^2+s^2(t)g$ with $s^{-1}\in L^1([0,\infty))$ which already contains the de-Sitter space and the power-law inflation metric. 
More generally, it is believed that generic future geodesically complete solutions of the Einstein equation with positive cosmological constant $\Lambda>0$ satisfy the assumptions of our theorems.
\end{Bsp}

Throughout this article we will employ the following notation:
We will use small Greek letters \(\alpha~\beta~\gamma\) for space-time indices,
small Latin indices \(i~j~k\) for spatial derivatives and
capital Latin indices \(I~J~K\) for indices on the Riemannian target.
We will denote spatial derivative by \(D\).
Moreover, we will make use of the usual summation convention,
that is we will always sum over repeated indices.

This article is organized as follows:
In the second section we introduce a suitable conformal transformation that we will be using to 
prove our main results. In the third section we establish the necessary energy estimate 
which is the key tool to prove the main result.

\section{Conformal Euler-Lagrange equations}
In the following we calculate how the energy functional \eqref{action-functional-quer} transforms under a certain conformal transformation of the metric on \(M\).
 
\begin{Lem}
If we transform the metric \(\tilde h=(Ns)^2h\) and the vector spinors via
\begin{align*}
\tilde\psi:=(Ns)^{\frac{1-n}{2}}\psi
\end{align*}
the energy functional \eqref{action-functional-quer} acquires the form
\begin{align}
\label{action-functional-h}
S(\phi,\psi,h)=\frac{1}{2}\int_{M}((Ns)^{n-2}|d\phi|^2+\langle\psi,i\D\psi\rangle-(Ns)^{2-n}\frac{1}{6}\langle\psi,R^P(\psi,\psi)\psi\rangle)\dv_{h}.
\end{align}
\end{Lem}
\begin{proof}
Under a conformal change of the metric \(\bar h=(Ns)^2h\) the volume elements transform as 
\(\dv_{\bar h}=(Ns)^n\dv_h\).
The Dirac operators transform as
\begin{align*}
\bar\D((Ns)^{-\frac{n-1}{2}}\bar\psi)=(Ns)^{-\frac{n+1}{2}}\overline{\D\psi}.
\end{align*}
Note that the twisted Dirac operator \(\D\) transforms in the same way as the standard Dirac operator \(\p\)
since the twist bundle \(\phi^\ast TP\) does not depend on the metric on \(M\).

Inserting our choice for \(\tilde\psi\) we get
\begin{align*}
\langle\tilde\psi,i\tilde\D\tilde\psi\rangle=(Ns)^{-n}\langle\psi,i\D\psi\rangle,
\end{align*}
which completes the proof.
\end{proof}

For the sake of completeness we calculate the critical points of \eqref{action-functional-h}.

\begin{Lem}[Critical points]
The critical points of \eqref{action-functional-h} are given by
\begin{align}
\label{system-dirac-wave-map-original}
\Box_h\phi=&(n-2)(Ns)^{-1}\nabla_{\nabla (Ns)}\phi-\frac{1}{2}(Ns)^{2-n}h^{\alpha\beta}R^P(\psi,i\partial_\alpha\cdot\psi)d\phi(\partial_\beta) \\
\nonumber &+\frac{1}{12}(Ns)^{4-2n}\langle(\nabla R^P)^\sharp(\psi,\psi)\psi,\psi\rangle, \\
\nonumber i\D\psi=&(Ns)^{2-n}\frac{1}{3}R^P(\psi,\psi)\psi.
\end{align}
Here, \(\Box:=-h^{\alpha\beta}\nabla_{\partial_\alpha}d\phi(\partial_\beta)\) denotes the wave map operator.
\end{Lem}
\begin{proof}
First, we vary the vector spinors \(\psi\) keeping the map \(\phi\) fixed.
More precisely, we consider a variation of \(\psi\colon(-\epsilon,\epsilon)\times M\to SM\otimes\phi^\ast TP\)
denoted by \(\psi_\lambda\) satisfying \(\frac{\nabla\psi_\lambda}{\partial\lambda}\big|_{\lambda=0}=\xi\).
We calculate
\begin{align*}
\frac{d}{d\lambda}\big|_{\lambda=0}\frac{1}{2}\int_M(&\langle\psi_\lambda,i\D\psi_\lambda\rangle-(Ns)^{2-n}\frac{1}{6}\langle\psi_\lambda,R^P(\psi_\lambda,\psi_\lambda)\psi_\lambda\rangle)\dv_{h}\\
=&\frac{1}{2}\int_M(\langle\xi,i\D\psi\rangle+\langle\psi,i\D\xi\rangle-(Ns)^{2-n}\frac{4}{6}\langle\xi,R^P(\psi,\psi)\psi\rangle)\dv_{h}\\
=&\int_M(\operatorname{Re}\langle\xi,i\D\psi\rangle-(Ns)^{2-n}\frac{1}{3}\langle\xi,R^P(\psi,\psi)\psi\rangle)\dv_{h},
\end{align*}
yielding the equation for the vector spinor.

Afterwards, we keep the vector spinors \(\psi\) fixed and consider a variation of \(\phi\colon(-\epsilon,\epsilon)\times M\to P\)
denoted by \(\phi_\lambda\) satisfying \(\frac{\partial\phi_\lambda}{\partial\lambda}\big|_{t=0}=\eta\).
We calculate
\begin{align*}
\frac{d}{d\lambda}\big|_{\lambda=0}\frac{1}{2}\int_M(Ns)^{n-2}|d\phi_\lambda|^2\dv_h=&\int_M(Ns)^{n-2}\langle\nabla\frac{\partial\phi_\lambda}{\partial\lambda},d\phi_\lambda\rangle \dv_h\big|_{\lambda=0} \\
=&\int_M\big((Ns)^{n-2}\langle\eta,\Box\phi\rangle
-(n-2)(Ns)^{n-3}\nabla(Ns)\langle\eta,\nabla\phi\rangle\big)\dv_h.
\end{align*}
Moreover, we have
\begin{align*}
\frac{d}{d\lambda}&\big|_{\lambda=0}\frac{1}{2}\int_M(\langle\psi,i\D\psi\rangle 
-(Ns)^{2-n}\frac{1}{6}\langle\psi,R^P(\psi,\psi)\psi\rangle)\dv_{h} \\
=&\frac{1}{2}\int_Mh^{\alpha\beta}\langle\psi,i R^{P}(\eta,d\phi(\partial_\alpha))\partial_\beta\cdot\psi\rangle \dv_h
+\frac{1}{12}\int_M(Ns)^{2-n}\big\langle\langle(\nabla R^P)^\sharp(\psi,\psi)\psi,\psi\rangle,\eta\big\rangle \dv_h,
\end{align*}
where we used the equation for \(\psi\). Finally, we have
\begin{align*}
h^{\alpha\beta}\langle\psi,i R^{P}(\eta,d\phi(\partial_\alpha))\partial_\beta\cdot\psi\rangle
=h^{\alpha\beta}\langle\eta, R^{P}(\psi,i\partial_\alpha\cdot\psi)d\phi(\partial_\beta)\rangle,
\end{align*}
which completes the proof.
\end{proof}

We want to turn the system \eqref{system-dirac-wave-map-original} into a system of two wave-type equations.
To this end we recall the following Weitzenböck formula for the twisted Dirac operator \(\D\).

\begin{Lem}
The square of the twisted Dirac operator \(\D\) satisfies
\begin{align}
\label{weitzenboeck}
\D^2=\Box+\frac{\scal^M}{4}+\frac{1}{2}h^{\alpha\gamma}h^{\beta\delta}\partial_\alpha\cdot\partial_\beta\cdot R^P(d\phi(\partial_\gamma),d\phi(\partial_\delta))\psi.
\end{align}
\end{Lem}
\begin{proof}
For a proof the we refer to \cite[Theorem 8.17]{MR1031992}.
\end{proof}

Making use of \eqref{weitzenboeck} we obtain the following rescaled Euler-Lagrange equations
\begin{align}
\label{wave-phi}\Box_h\phi=&(n-2)(Ns)^{-1}\nabla_{\nabla (Ns)}\phi-\frac{1}{2}(Ns)^{2-n}h^{\alpha\beta}R^P(\psi,i\partial_\alpha\cdot\psi)d\phi(\partial_\beta) \\
\nonumber &+\frac{1}{12}(Ns)^{4-2n}\langle(\nabla R^P)^\sharp(\psi,\psi)\psi,\psi\rangle, \\
\label{wave-psi}\Box_h\psi=&-\frac{\scal^M}{4}\psi-\frac{1}{2}h^{\alpha\gamma}h^{\beta\delta}\partial_\alpha\cdot\partial_\beta\cdot R^P(d\phi(\partial_\gamma),d\phi(\partial_\delta))\psi
+\frac{i}{3}\nabla(R^P(\psi,\psi)(Ns)^{2-n})\cdot\psi \\
&\nonumber+\frac{1}{9}(Ns)^{4-2n}R^P(\psi,\psi)R^P(\psi,\psi)\psi,
\end{align}
where \(\Box_h\) denotes the corresponding wave-operator.

In terms of local coordinates this system acquires the form (\(I=1,\ldots,\dim P\))
\begin{align*}
s^2\partial^2_t\phi^I&+s\dot s\partial_t\phi^I+\frac{1}{2}s^2\tr\dot g\partial_t\phi^I+D^\ast D\phi^I \\
=&-(n-2)\dot ss\partial_t\phi^I-(n-2)s^2N^{-1}\partial_tN\partial_t\phi^I 
+(n-2)N^{-1}D_{\partial_i}ND_{\partial_i}\phi^I \\
&-\frac{1}{2}(Ns)^{2-n}R^I_{JKL}h^{\alpha\beta}\langle\psi^K,i\partial_\alpha\cdot\psi^L\rangle\partial_\beta\phi^J \\
&+\frac{1}{12}(Ns)^{4-2n}G^{IJ}\nabla_JR_{KLMN}\langle\psi^K,\psi^M\rangle\langle\psi^L,\psi^N\rangle, \\
s^2\nabla^2_t\psi^I&+s\dot s\nabla_t\psi^I+\frac{1}{2}s^2\tr  \dot g\nabla_t\psi^I+D^\ast D\psi^I \\
=&-\frac{\scal^M}{4}\psi^I
-\frac{1}{2}h^{\alpha\gamma}h^{\beta\delta}R^I_{JKL}\partial_\gamma\phi^K\partial_\delta\phi^L\partial_\alpha\cdot\partial_\beta\cdot\psi^J\\
&+\frac{i}{3}h^{\alpha\beta}\nabla_{\partial_\alpha}\big((Ns)^{2-n}R^I_{JKL}\langle\psi^J,\psi^L\rangle\partial_{y^I}\big)\partial_\beta\cdot\psi^K\\
&+\frac{1}{9}(Ns)^{4-2n}R^I_{JKL}R^K_{MRS}\langle\psi^J,\psi^L\rangle\langle\psi^M,\psi^S\rangle\psi^R.
\end{align*}

\begin{Bem}
Note that the system \eqref{phi-dwmap-quer}, \eqref{psi-dwmap-quer} on the manifold
\((M,\tilde h)\) is equivalent to the system
\eqref{system-dirac-wave-map-original} on the conformally transformed manifold \((M,h)\).
From now on we will use the system \eqref{wave-phi}, \eqref{wave-psi} which follows from
\eqref{system-dirac-wave-map-original}.
\end{Bem}

\section{The energy estimate}
In this section we first develop several formulas that are useful in the study
of energy estimates for sections in arbitrary vector bundles. Later on, we will
apply these techniques to the cases of wave maps and Dirac-wave maps with curvature term.

\subsection{Energies on general vector bundles}
Let $\Sigma^n$ be a manifold, $s:\R\to\R_+$ smooth and $g_t$, $t\in\R$ be a smooth family of Riemannian metrics on $\Sigma$. We consider the (globally hyperbolic) Lorentzian manifold 
\begin{align*}
(M,h)=(\R\times\Sigma, -s(t)^{-2}dt^2+g_t).
\end{align*}
The non-vanishing Christoffel symbols of this metric are given by
\begin{align*}
\Gamma(h)_{00}^0=-\frac{\dot{s}}{s},\qquad\Gamma(h)_{ij}^0=\frac{1}{2}s^2\dot{g}_{ij},\qquad \Gamma(h)_{i0}^j=\Gamma(h)_{0i}^j=\frac{1}{2}g^{jk}\dot{g}_{ik},\qquad \Gamma(h)_{ij}^k=\Gamma(g)_{ij}^k.
\end{align*}
Let $V$ be a Riemannian vector bundle over $M$ which is equipped with a metric connection $\nabla$. As usual, iterating yields a map 
\begin{align*}
\nabla^k:\Gamma(V)\to \Gamma(T^*M^{\otimes k}\otimes V).
\end{align*} 
We define the associated wave operator $\Box:\Gamma(V)\to\Gamma(V)$ by the sign convention such that
\begin{align*}
\Box\xi=-h^{\alpha\beta}\nabla^2_{\alpha\beta}\xi
\end{align*}
for \(\xi\in\Gamma(V)\).
The covariant derivative $\nabla$ restricts for each $t\in\R$ to the spatial covariant derivative $D$ which yields a map
\begin{align*}
D:\Gamma(V)\to \Gamma(\pi^*(T^*\Sigma)\otimes V),
\end{align*}
where $\pi:M\to\Sigma$ is the canonical projection. In the following, we will write $T^*\Sigma$ instead of $\pi^*(T^*\Sigma)$ for notational convenience. The covariant derivative naturally extends as a map
 \begin{align*}D:\Gamma(T^*\Sigma^{\otimes k}\otimes V)\to \Gamma(T^*\Sigma^{\otimes k+1}\otimes V)
 \end{align*}
 by defining
 \begin{align*}
 (D_X\xi)(X_1,\ldots,X_k)={}^{V}\nabla_X(\xi(X_1,\ldots,X_k))-\sum_{i=1}^k\xi(X_1,\ldots,D_XX_i,\ldots,X_k),
 \end{align*}
 where for each $t\in\R$, $D_XX_i$ denotes the covariant derivative of $g_t$.
Furthermore, we define its formal adjoint $D^*:\Gamma(T^*\Sigma^{\otimes k+1}\otimes V)\to \Gamma(T^*\Sigma^{\otimes k}\otimes V)$ by
\begin{align*}
D^*\xi(X_1,\ldots,X_k)=-g^{ij}D_{\partial_i}\xi(\partial_j,X_1,\ldots,X_k).
\end{align*}
Finally, we define a covariant time-derivative $\nabla_t:\Gamma(T^*\Sigma^{\otimes k}\otimes V)\to\Gamma(T^*\Sigma^{\otimes k}\otimes V)$
by
 \begin{align*}
 (\nabla_t\xi)(X_1,\ldots,X_k)=\nabla^V_{\partial_t}(\xi(X_1,\ldots,X_k))-\sum_{i=1}^k\xi(X_1,\ldots,\nabla_{\partial_t}X_i,\ldots,X_k).
 \end{align*}
Observe that this definition makes sense as $\nabla_{\partial_t}X_i\in\Gamma(TM)$ is always tangential to $\Sigma$ due to the structure of the metric $h$.
A quick computation shows that the wave operator can be written as
\begin{align}\label{waveoperator}
\Box\xi=s^2\nabla_t\nabla_t\xi+\dot{s}s\nabla_t\xi+\frac{1}{2}s^2\trace_g\dot{g}\cdot\nabla_t\xi+D^*D\xi.
\end{align}

\begin{Lem}\label{productrule}We use the notations from above.
Let $\langle.,.\rangle$ be the natural $t$-dependent scalar product induced on $\Gamma(T^*\Sigma^{\otimes k}\otimes V)$. Let $\xi,\eta\in\Gamma(V)$. Then we have the product rule
\begin{align*}
\partial_t\langle D^k\xi,D^k\eta\rangle=
\langle\nabla_t D^k\xi,D^k\eta\rangle+\langle D^k\xi,\nabla_t D^k\eta\rangle.
\end{align*}
\end{Lem}
\begin{proof}
Let $t_0$ be a fixed time and let $\left\{e_1,\ldots,e_n\right\}$ be an orthonormal basis of $(\Sigma,g_{t_0})$ at $x\in \Sigma$. Then the scalar product can also be written as
\begin{align*}
\langle D^k\xi,D^k\eta\rangle(t_0,x)=\sum_{i_1,\ldots,i_k=1}^n \langle D^k_{e_{i_1},\ldots,e_{i_k}}\xi,D^k_{e_{i_1},\ldots,e_{i_k}}\eta\rangle(t_0,x),
\end{align*}
where the scalar products on the right hand side are the ones on $V$. Now think of $\left\{e_1,\ldots,e_n\right\}$ as an orthonormal system in $T_{(t_0,x)}M$. Parallel transport along the curve $t\mapsto (t,x)$ yields orthonormal systems $\left\{e_1,\ldots,e_n\right\}$ for each $T_{(t,x)}M$. Because $\nabla_{\partial_t}\partial_t=-\frac{\dot{s}}{s}\partial_t$, $h(\partial_t,e_i)\equiv0$ for each $i\in\left\{1,\ldots,n\right\}$ so that we in fact obtain orthonormal bases of $T_x\Sigma$ with respect to the metric $g_t$ for each $t\in \R$. Therefore we get
\begin{align*}
\partial_t\langle D^k\xi,D^k\eta\rangle&=\sum_{i_1,\ldots,i_k=1}^n \langle\nabla_t (D^k_{e_{i_1},\ldots,e_{i_k}}\xi),D^k_{e_{i_1},\ldots,e_{i_k}}\eta\rangle+\sum_{i_1,\ldots,i_k=1}^n \langle D^k_{e_{i_1},\ldots,e_{i_k}}\xi,\nabla_t (D^k_{e_{i_1},\ldots,e_{i_k}}\eta)\rangle\\
&=\sum_{i_1,\ldots,i_k=1}^n \langle\nabla_t D^k_{e_{i_1},\ldots,e_{i_k}}\xi,D^k_{e_{i_1},\ldots,e_{i_k}}\eta\rangle+\sum_{i_1,\ldots,i_k=1}^n \langle D^k_{e_{i_1},\ldots,e_{i_k}}\xi,\nabla_t D^k_{e_{i_1},\ldots,e_{i_k}}\eta\rangle\\
&=\langle\nabla_t D^k\xi,D^k\eta\rangle+\langle D^k\xi,\nabla_t D^k\eta\rangle.
\end{align*}
\end{proof}
\begin{Bem}
In the above lemma, the special structure of the metric $h$ is essential.
\end{Bem}
For a section $\xi\in \Gamma(V)$, we now define the $k$th energy density
\begin{align*}
e_k(\xi)=s^2|D^k\nabla_t\xi|^2+|D^{k+1}\xi|^2,
\end{align*}
which is a nonnegative function on $M$ and the $k$th energy as
\begin{align*}
E_k(\xi)(t)=\int_{\Sigma}e_k(\xi)\dv_{g_t},
\end{align*}
where $\dv_{g_t}$ is the volume element of the metric $g_t$ and the integral has to be understood as an integral over $\left\{t\right\}\times\Sigma$.
\begin{Prop}\label{energy-evolution-general}
Let $\xi\in\Gamma(V)$ be spacelike compactly supported. Then its $k$th energy satisfies the evolution equation
\begin{align*}
\frac{d}{dt} E_k(\xi)&=2\int_{\Sigma}\langle D^k\Box\xi,D^k\nabla_t\xi\rangle\dv_{g_t}+\frac{1}{2}\int_{\Sigma}(|D^{k+1}\xi|^2-s^2|D^k\nabla_t\xi|^2)\trace_g\dot{g}\dv_{g_t}\\
&\qquad+2s^2\int_{\Sigma}\langle D^k\nabla_t\xi,[\nabla_t,D^k]\nabla_t\xi\rangle\dv_{g_t}+2\int_{\Sigma}\langle D^k\nabla_t\xi,[D^*D,D^k]\xi\rangle\dv_{g_t}\\
&\qquad+2\int_{\Sigma}\langle D^{k+1}\xi,[\nabla_t,D^{k+1}]\xi\rangle\dv_{g_t}.
\end{align*}
\end{Prop}
\begin{proof}
We compute
\begin{align*}
\frac{d}{dt}E_k(\xi)&=2\dot{s}s\int_{\Sigma}|D^k\nabla_t\xi|^2\dv_{g_t}+\frac{1}{2}\int_{\Sigma}e_k(\xi)\trace_g\dot{g}\dv_{g_t}\\ 
&\qquad+2s^2\int_{\Sigma}\langle\nabla_tD^k\nabla_t\xi,D^k\nabla_t\xi\rangle\dv_{g_t}+2\int_{\Sigma}\langle\nabla_tD^{k+1}\xi,D^{k+1}\xi\rangle\dv_{g_t}\\
&=2\dot{s}s\int_{\Sigma}|D^k\nabla_t\xi|^2\dv_{g_t}+\frac{1}{2}\int_{\Sigma}e_k(\xi)\trace_g\dot{g}\dv_{g_t}\\ 
&\qquad+2s^2\int_{\Sigma}\langle D^k\nabla_t\nabla_t\xi,D^k\nabla_t\xi\rangle\dv_{g_t}+2s^2\int_{\Sigma}\langle[\nabla_t,D^k]\nabla_t\xi,D^k\nabla_t\xi\rangle\dv_{g_t}\\
&\qquad+2\int_{\Sigma}\langle D^{k+1}\nabla_t\xi,D^{k+1}\xi\rangle\dv_{g_t}+2\int_{\Sigma}\langle[\nabla_t,D^{k+1}]\xi,D^{k+1}\xi\rangle\dv_{g_t}.
\end{align*}
Moreover, we have
\begin{align*}
2\int_{\Sigma}\langle D^{k+1}\nabla_t\xi,D^{k+1}\xi\rangle\dv_{g_t}&=2\int_{\Sigma}\langle D^k\nabla_t\xi,D^*D^{k+1}\xi\rangle\dv_{g_t}\\
&= 2\int_{\Sigma}\langle D^k\nabla_t\xi,D^k D^*D\xi\rangle\dv_{g_t}+2\int_{\Sigma}\langle D^k\nabla_t\xi,[D^{*}D,D^k]\xi\rangle\dv_{g_t}.
\end{align*}
The statement follows by combining both equalities and using \eqref{waveoperator}.
\end{proof}

In order to derive energy estimates for the vector spinors we have to take into account the
special structure of the scalar product of the spinor bundle \(SM\). The natural geometric scalar product
which is invariant under the spin group is not positive definite and thus not a good candidate
for analytic purposes. To obtain a positive definite scalar product one has to Clifford multiply
the second factor with the timelike unit vector field \(e_0\), see \cite{MR701244}, which in our setup is given by \(e_0=s\partial_t\).
Using our geometric setup we find
\begin{align*}
\nabla_te_0=0,\qquad \nabla_ie_0=\frac{1}{2}sg^{jk}\dot g_{ik}\partial_{x^j}.
\end{align*}
Hence, we define the following \(k\)th energy density for the vector spinors \(\psi\)
\begin{align*}
e_k(\psi)=s^2\langle D^k\nabla_t\psi,e_0\cdot D^k\nabla_t\psi\rangle+\langle D^{k+1}\psi,e_0\cdot D^{k+1}\psi\rangle,
\end{align*}
which is a nonnegative function. In the following we will always employ the positive definite scalar product
without mentioning it explicitly. However, since \(e_0\) is not parallel with respect to the spatial coordinates
the \(k\)th energy of the spinor satisfies
\begin{align}
\label{kenergy-spinor}
\frac{d}{dt} E_k(\psi)&=2\int_{\Sigma}\langle D^k\Box\psi, D^k\nabla_t\psi\rangle\dv_{g_t}+\frac{1}{2}\int_{\Sigma}(|D^{k+1}\psi|^2-s^2|D^k\nabla_t\psi|^2)\trace_g\dot{g}\dv_{g_t}\\
\nonumber&\qquad+2s^2\int_{\Sigma}\langle D^k\nabla_t\psi,[\nabla_t,D^k]\nabla_t\psi\rangle\dv_{g_t}+2\int_{\Sigma}\langle D^k\nabla_t\psi,[D^*D,D^k]\psi\rangle\dv_{g_t}\\
\nonumber&\qquad+2\int_{\Sigma}\langle D^{k+1}\psi,[\nabla_t,D^{k+1}]\psi\rangle\dv_{g_t}
+2\int_{\Sigma}\langle D^{k}\psi,e_0\cdot(D^\ast e_0)\cdot D^{k+1}\psi\rangle\dv_{g_t}.
\end{align}

\begin{Lem}
Let $\xi\in\Gamma(T^*\Sigma^{\otimes k}\otimes V)$. Then we have the identity
\begin{align*}
([\nabla_t,D]\xi)(X,X_1,\ldots,X_k)&=R^V_{\partial_t,X}(\xi(X_1,\ldots,X_k))-\frac{1}{2}D_{\dot{g}(X)}\xi(X_1,\ldots,X_k)\\&\qquad+\frac{1}{2}\sum_{i=1}^k\xi(X_1,\ldots,D_X\dot{g}(X_i),\ldots,X_k).
\end{align*}
\end{Lem}
\begin{proof}
At first, we compute
\begin{align*}
&(\nabla_tD\xi)(X,X_1,\ldots,X_k)=\nabla^V_{\partial_t}(D\xi(X,X_1,\ldots,X_k))-(D\xi)(\nabla_{\partial_t}X,X_1,\ldots,X_k)\\
&\qquad-\sum_{i=1}^k(D\xi)(X,X_1,\ldots,\nabla_{\partial_t}X_i,\ldots,X_k)\\
&=\nabla^V_{\partial_t}(\nabla^V_X(\xi(X_1,\ldots,X_k))-
\sum_{i=1}^k\xi(X_1,\ldots,D_XX_i,\ldots X_k))\\
&\qquad -(D_{\nabla_{\partial_t}X}\xi)(X_1,\ldots,X_k)
-\sum_{i=1}^k(D_X\xi)(X_1,\ldots,\nabla_{\partial_t}X_i,\ldots,X_k)\\
&=\nabla^V_{\partial_t}(\nabla^V_X(\xi(X_1,\ldots,X_k)))
-\sum_{i=1}^k(\nabla_t\xi)(X_1,\ldots,D_XX_i,\ldots X_k)\\
&\qquad-\sum_{i=1}^k\xi(X_1,\ldots,\nabla_{\partial_t}D_XX_i,\ldots X_k)
-\sum_{\substack{i,j=1\\i\neq j}}^n
\xi(X_1,\ldots,D_XX_i,\ldots,\nabla_{\partial_t}X_j,\ldots X_k)\\&\qquad  -(D_{\nabla_{\partial_t}X}\xi)(X_1,\ldots,X_k)
-\sum_{i=1}^k(D_X\xi)(X_1,\ldots,\nabla_{\partial_t}X_i,\ldots,X_k).
\end{align*}
Similarly, we find
\begin{align*}
&(D(\nabla_t\xi))(X,X_1,\ldots,X_k)=\nabla^V_X((\nabla_t\xi)(X_1,\ldots,X_k))-\sum_{i=1}^k(\nabla_t\xi)(X_1,\ldots,D_XX_i,\ldots,X_k)\\
&=\nabla^V_X(\nabla^V_{\partial_t}(\xi(X_1,\ldots,X_k))-\sum_{i=1}^k\xi(X_1,\ldots,\nabla_{\partial_t}X_i,\ldots,X_k))\\&\qquad-\sum_{i=1}^k(\nabla_t\xi)(X_1,\ldots,D_XX_i,\ldots,X_k)\\
&=\nabla^V_X(\nabla^V_{\partial_t}(\xi(X_1,\ldots,X_k)))
-\sum_{i=1}^k(D_X\xi)(X_1,\ldots,\nabla_{\partial_t}X_i,\ldots X_k)\\
&\qquad-\sum_{i=1}^k\xi(X_1,\ldots,D_X\nabla_{\partial_t}X_i,\ldots X_k)
-\sum_{\substack{i,j=1\\i\neq j}}^n
\xi(X_1,\ldots,D_XX_i,\ldots,\nabla_{\partial_t}X_j,\ldots X_k)\\&\qquad-\sum_{i=1}^k(\nabla_t\xi)(X_1,\ldots,D_XX_i,\ldots,X_k).
\end{align*}
Summing up, we obtain
\begin{align*}
([\nabla_t,D]\xi)(X,X_1,\ldots,X_k)&=\nabla^V_{\partial_t}(\nabla^V_X(\xi(X_1,\ldots,X_k)))
-\nabla^V_X(\nabla^V_{\partial_t}(\xi(X_1,\ldots,X_k)))\\
&\qquad-\sum_{i=1}^k\xi(X_1,\ldots,\nabla_{\partial_t}D_XX_i-D_X\nabla_{\partial_t}X_i,\ldots X_k)\\&\qquad-(D_{\nabla_{\partial_t}X}\xi)(X_1,\ldots,X_k).
\end{align*}
As this expression is tensorial in $X$ and the $X_i$, we may assume that the components of these vector fields with respect to a chart of $\Sigma$ are independent of time.
By raising an index with respect to $g_t$ we can think of $\dot{g}_t$ as an endomorphism on $T^*\Sigma$.
Then the formulas for the Christoffel symbols imply $\nabla_{\partial_t}X=\nabla_X\partial_t=\frac{1}{2}\dot{g}(X)$ and the same holds for $X_i$. 
Moreover, we have $[\partial_t,X]=[\partial_t,X_i]=0$ and 
\begin{align*}
\nabla_{\partial_t}D_XX_i-D_X\nabla_{\partial_t}X_i=\frac{1}{2}\big(\dot{g}(D_XX_i)-D_X(\dot{g}(X_i))\big)=-\frac{1}{2}D_X\dot{g}(X_i).
\end{align*}
By putting these facts together, we immediately get the statement of the Lemma.
\end{proof}

In the following we will often make use of the so-called \(\star\) notation.
More precisely, we will use a \(\star\) to denote various contractions
between the objects involved.

We apply the general formula from above
in the case where we have a map $\phi\in C^{\infty}(M,P)$,
$[\nabla_t,D]\phi\in \Gamma(T^*\Sigma\otimes \phi^\ast TP)$ and
we get
\begin{align*}
([\nabla_t,D]\phi)(X)=-\frac{1}{2}d\phi(\dot{g}(X)).
\end{align*}
More generally, if $E=\phi^*TP$ in the Lemma above, we get
\begin{align*}
([\nabla_t,D]\xi)(X,X_1,\ldots,X_k)&=R^P(d\phi(\partial_t),d\phi(X))(\xi(X_1,\ldots,X_k))-\frac{1}{2}D_{\dot{g}(X)}\xi(X_1,\ldots,X_k)\\&\qquad+\frac{1}{2}\sum_{i=1}^k\xi(X_1,\ldots,D_X\dot{g}(X_i),\ldots,X_k)\\
&=R^P(\nabla_t\phi,D_X\phi)(\xi(X_1,\ldots,X_k))-\frac{1}{2}D_{\dot{g}(X)}\xi(X_1,\ldots,X_k)\\&\qquad+\frac{1}{2}\sum_{i=1}^k\xi(X_1,\ldots,D_X\dot{g}(X_i),\ldots,X_k).
\end{align*}
By an iteration argument, we find
\begin{align*}
[\nabla_t,D^k]\nabla_t\phi&=\sum_{\sum l_i+\sum {m_j}=k-1}{}^{G}\nabla^{l_1}R^P\star\underbrace{D^{m_1+1}\phi\star\ldots\star D^{m_{l_1}+1}\phi}_{l_1-\text{times}}\star D^{l_2}\nabla_t\phi\star D^{l_3+1}\phi\star D^{l_4}\nabla_t\phi\\
&\qquad+\sum_{l=0}^{k-1}D^{k-l}\dot{g}\star D^{l}\nabla_t\phi,\\
[\nabla_t,D^{k+1}]\phi&=\sum_{\sum l_i+\sum {m_j}=k-1}{}^{G}\nabla^{l_1}R^N\star\underbrace{D^{m_1+1}\phi\star\ldots\star D^{m_{l_1}+1}\phi}_{l_1-\text{times}}\star D^{l_2}\nabla_t\phi\star D^{l_3+1}\phi\star D^{l_4+1}\phi\\
&\qquad+\sum_{l=0}^{k}D^{k-l}\dot{g}\star D^{l+1}\phi.
\end{align*}
In the case of a vector spinor $\psi\in\Gamma(SM\otimes \phi^*TP)$, the formulas are
\begin{align*}
([\nabla_t,D]\psi)(X)=R^{SM}({\partial_t,X})\cdot\psi+R^P(\nabla_t\phi,D_X\phi)\psi-\frac{1}{2}D\phi(\dot{g}(X))
\end{align*}
and
\begin{align*}
([\nabla_t,D]\xi)(X,X_1,\ldots,X_k)&=R^{SM}(\partial_t,X)\cdot\xi+R^P(\nabla_t\phi,D_X\phi)(\xi(X_1,\ldots,X_k))\\&\qquad
-\frac{1}{2}D_{\dot{g}(X)}\xi(X_1,\ldots,X_k)+\frac{1}{2}\sum_{i=1}^k\xi(X_1,\ldots,D_X\dot{g}(X_i),\ldots,X_k).
\end{align*}

By iterating the formula from above, we get 
\begin{align*}
[\nabla_t,D^k]\nabla_t\psi&=
\sum_{l=0}^{k-1}D^l(R^{SM}(\partial_t,.))\star D^{k-1-l}\nabla_t\psi
+\sum_{l=0}^{k-1}D^{k-l}\dot{g}\star D^{l+1}\nabla_t\psi\\
&\quad+\sum_{\sum l_i+\sum {m_j}=k-1}{}^{G}\nabla^{l_1}R^P\star\underbrace{D^{m_1+1}\phi\star\ldots\star D^{m_{l_1}+1}\phi}_{l_1-\text{times}}\star D^{l_2}\nabla_t\phi\star D^{l_3+1}\phi\star D^{l_4}\nabla_t\psi
\end{align*}
and
\begin{align*}
[\nabla_t,D^{k+1}]\psi&=\sum_{l=0}^{k}D^l(R^{SM}(\partial_t,.))\star D^{k-l}\psi+\sum_{l=0}^{k}D^{k-l+1}\dot{g}\star D^{l+1}\psi\\
&\quad	+\sum_{\sum l_i+\sum {m_j}=k}{}^{G}\nabla^{l_1}R^P\star\underbrace{D^{m_1+1}\phi\star\ldots\star D^{m_{l_1}+1}\phi}_{l_1-\text{times}}\star D^{l_2}\nabla_t\phi\star D^{l_3+1}\phi\star D^{l_4}\psi.
\end{align*}
\begin{Lem}
Let $\xi\in\Gamma(T^*\Sigma^{\otimes k}\otimes V)$. Then we have the identity
\begin{align*}
([D^*D,D]\xi)(X,X_1,\ldots,X_k)&=-D_{\ric(X)}\xi(X_1,\ldots,X_k)
-2\sum_{j=1}^kD_{e_i}\xi(X_1,\ldots R_{X,e_i}X_j,\ldots,X_k)\\
&\qquad+ D_{e_i}R^V_{X,e_i}(\xi(X_1,\ldots,X_k))
+2R^V_{X,e_i}(D_{e_i}\xi(X_1,\ldots,X_k))\\
&\qquad-\sum_{j=1}^k\xi(X_1,\ldots D_{e_i}R_{X,e_i}X_j,\ldots,X_k),
\end{align*}
where $\left\{e_i\right\}_{1\leq i\leq n-1}$ is a local orthonormal frame. Here and throughout the proof, we sum over $i$.
\end{Lem}
\begin{proof}
A direct computation yields
\begin{align*}
([D^*D,D]\xi)&(X,X_1,\ldots,X_k)=D^3_{X,e_i,e_i}\xi(X_1,\ldots,X_k)-D^3_{e_i,e_i,X}\xi(X_1,\ldots,X_k)\\&=D^3_{X,e_i,e_i}\xi(X_1,\ldots,X_k)-D^3_{e_i,X,e_i}\xi(X_1,\ldots,X_k)\\&\qquad+D^3_{e_i,X,e_i}\xi(X_1,\ldots,X_k)-D^3_{e_i,e_i,X}\xi(X_1,\ldots,X_k)\\
&=R_{X,e_i}D\xi(e_i,X_1,\ldots,X_k)+D_{e_i}R_{X,e_i}\xi(X_1,\ldots,X_k)\\
&=-D_{\ric(X)}\xi(X_1,\ldots,X_k)
-\sum_{j=1}^kD_{e_i}\xi(X_1,\ldots R_{X,e_i}X_j,\ldots,X_k)\\
&\qquad+R_{X,e_i}(D_{e_i}\xi(X_1,\ldots,X_k))\\
&\qquad+D_{e_i} (R^V\circ\xi)(X,e_i,X_1,\ldots,X_k)
-D_{e_i}\sum_{j=1}^k\xi(X_1,\ldots R_{X,e_i}X_j,\ldots,X_k)\\
&=-D_{\ric(X)}\xi(X_1,\ldots,X_k)
-2\sum_{j=1}^kD_{e_i}\xi(X_1,\ldots R_{X,e_i}X_j,\ldots,X_k)\\
&\qquad+ D_{e_i}R^V_{X,e_i}(\xi(X_1,\ldots,X_k))
+2R_{X,e_i}(D_{e_i}\xi(X_1,\ldots,X_k))\\
&\qquad-\sum_{j=1}^k\xi(X_1,\ldots D_{e_i}R_{X,e_i}X_j,\ldots,X_k).
\end{align*}
Here, $R^V\circ\xi\in \Gamma(T^*\Sigma^{\otimes k+2}\otimes V)$ is defined as 
\begin{align*}
(R^V\circ\xi)(X,Y,X_1,\ldots,X_k)=R_{X,Y}^V(\xi(X_1,\ldots,X_k))
\end{align*}
and $DR^{V}$ is the covariant derivative of the curvature endomorphism on $V$ restricted to vectors tangential to $\Sigma$. In other words,
\begin{align*}
D_{X}R^V_{Y,Z}\xi:=D_X(R^V_{Y,Z}\xi)-R^V_{D_XY,Z}\xi-R^V_{Y,D_X,Z}\xi-R^V_{Y,Z}(\nabla_X\xi).\\
\qedhere
\end{align*}
\end{proof}
In the  case $\phi\in C^{\infty}(M,P)$,
we obtain
\begin{align*}
([D^*D,D]\phi)(X)=-D\phi(\ric^M(X))+R^P(D\phi(X),D\phi(e_i))(D\phi(e_i))
\end{align*}
and for $\xi\in \Gamma(T^*\Sigma^{\otimes k}\otimes \phi^*TP)$, we have
\begin{align*}
([&D^*D,D]\xi)(X,X_1,\ldots X_k)=2\sum_{l=1}^kD_{e_i}\xi(X_1,\ldots, R^M(e_i,X)X_l,\ldots X_k)\\&\qquad+\sum_{l=1}^k\xi(X_1,\ldots ,D_{e_i}R^M(e_i,X)X_l,\ldots X_k)-D_{\ric(X)}\xi(X_1,\ldots,X_k)\\
&\qquad +D_{D\phi(e_i)}R^P(D\phi(X),D\phi(e_i))(\xi(X_1,\ldots X_k))
+R^P(D^2\phi({e_i,X}),D\phi(e_i))(\xi( X_1,\ldots X_k))\\
&\qquad +R^P(D\phi(X),D^2\phi(e_i,e_i))(\xi(X_1,\ldots X_k))+2R^P(D\phi(X),D\phi(e_i))(D_{e_i}\xi(X_1,\ldots,X_k)).
\end{align*}
Thus by iteration we get
\begin{align*}
[D^*D,D^k]\phi&=\sum_{l=0}^{k-1}D^{l}R^M\star D^{k-l}\phi\\&\quad+\sum_{\sum l_i+\sum {m_j}=k-1}{}^{G}\nabla^{l_1}R^P\star\underbrace{D^{m_1+1}\phi\star\ldots\star D^{m_{l_1}+1}\phi}_{l_1-\text{times}}
\star D^{l_2+1}\phi\star D^{l_3+1}\phi\star D^{l_4+1}\phi.
\end{align*}
In the case of vector spinors, that is $\psi\in \Gamma(SM\otimes \phi^*TP)$, we obtain
\begin{align*}
([D^*D,D]\psi)(X)&=-D_{\ric(X)}\psi+2R^P(D\phi(X),D\phi(e_i))(D\psi(e_i))+2R^{SM}_{X,e_i}\cdot D_{e_i}\psi+D_{e_i}R^{SM}_{X,e_i}\cdot \psi\\
&\qquad+\nabla_{D\phi(e_i)}R^P(D\phi(X),D\phi(e_i))\psi+
R^P(D^2\phi(e_i,X),D\phi(e_i))\psi\\
&\qquad+R^P(D\phi(X),D^2\phi(e_i,e_i))\psi
\end{align*}
and for $\xi\in \Gamma(T^*\Sigma^{\otimes k}\otimes SM\otimes \phi^*TP)$, we have
\begin{align*}
([D^*D,D]\xi)(X,X_1,\ldots,X_k)&=-D_{\ric(X)}\xi(X_1,\ldots,X_k)\\ &\qquad+2R^P(D\phi(X),D\phi(e_i))(D_{e_i}\psi(X_1,\ldots,X_k))\\& \qquad+2R^{SM}_{X,e_i}\cdot D_{e_i}\psi(X_1,\ldots,X_k)+D_{e_i}R^{SM}_{X,e_i}\cdot \psi(X_1,\ldots,X_k)\\
&\qquad+\nabla_{D\phi(e_i)}R^P(D\phi(X),D\phi(e_i))(\psi(X_1,\ldots,X_k))\\&\qquad+
R^P(D^2\phi(e_i,X),D\phi(e_i))(\psi(X_1,\ldots,X_k))\\
&\qquad+R^P(D\phi(X),D^2\phi(e_i,e_i))(\psi(X_1,\ldots,X_k))\\
&\qquad+2\sum_{l=1}^kD_{e_i}\xi(X_1,\ldots, R^M(e_i,X)X_l,\ldots X_k)\\&\qquad+\sum_{l=1}^k\xi(X_1,\ldots ,D_{e_i}R^M(e_i,X)X_l,\ldots X_k).
\end{align*}
By iteration, we then obtain
\begin{align*}
[D^*D,D^k]\psi=&\sum_{l=0}^{k-1}D^{l}R^M\star D^{k-l}\psi+\sum_{l=0}^{k}D^{l}R^{SM}\star D^{k-l}\psi
\\&\quad+\sum_{\sum l_i+\sum {m_j}=k}{}^{G}\nabla^{l_1}R^N\star\underbrace{D^{m_1+1}\phi\star\ldots\star D^{m_{l_1}+1}\phi}_{l_1-\text{times}}\star D^{l_2+1}\phi\star D^{l_3+1}\phi\star D^{l_4}\psi.
\end{align*}

We conclude this section with a very important lemma,
which we will frequently make use of when deriving energy estimates.

\begin{Lem}\label{sobolev_multiplication}
Let $(M^m,g)$ be an $m$-dimensional Riemannian manifold, $r$ a natural number satisfying $r>\frac{m}{2}$, $E\to M$ a Riemannian vector bundle with a connection, $k\in \N$, $\xi_1,\ldots \xi_k\in H^r(E)$ 
and $l_i$ $i=1,\ldots k$ natural numbers satisfying $\sum_i l_i\leq r$. Then the following inequality holds
\begin{align*}
\int_M \sff_{i=1}^k|\nabla^{l_i}\xi_i|^2\dv\leq C_{Sob}\cdot\sff_{i=1}^k\left\| \xi_i\right\|_{H^r}.
\end{align*}
If in addition, we have $\xi_{k+1}\in H^{r-1}(E)$ and  $\sum_{i=1}^{k+1} l_i\leq r-1$, then 
\begin{align*}
	\int_M \sff_{i=1}^{k+1}|\nabla^{l_i}\xi_i|^2\dv\leq C_{Sob}\cdot\sff_{i=1}^k\left\| \xi_i\right\|_{H^r}\cdot\left\| \xi_{k+1}\right\|_{H^{r-1}} .
\end{align*}
Here, \(C_{sob}=C_{sob}(g,l_1,\ldots,l_k,r)\) depends on the constant from the Sobolev embedding on \(M\) and
the numbers \(l_1,\ldots,l_k,r\).
\end{Lem}

\begin{proof}We prove the first inequality, the second one is shown very similarly.
Choose $p_i\in [1,\infty]$, $i=1,\ldots,l_k$ such that $\frac{1}{p_i}> \frac{1}{2}-\frac{r-l_i}{n}$ so that $\left\| \nabla^{l_i}\xi_i\right\|_{L^{p_i}}\leq C\left\| \xi_i\right\|_{H^r}$ by Sobolev embedding. 
Furthermore,
\begin{align*}
\sum_{i=1}^k\left(\frac{1}{2}-\frac{r-l_i}{m}\right)
=\frac{k}{2}-\frac{kr}{m}+\frac{1}{m}(\sum_{i=1}^kl_i)=\frac{1}{2}+(k-1)\left(\frac{1}{2}-\frac{r}{m}\right)+\frac{1}{m}\left(\sum_{i=1}^kl_i-r\right)<\frac{1}{2}
\end{align*}
due to our assumptions. Therefore, we can choose the $p_i$ such that $\sum_{i=1}^k\frac{1}{p_i}=\frac{1}{2}$. An application of the H\"{o}lder inequality finishes the proof of the inequality.
\end{proof}

\subsection{Energy of the map}
We define the $k$th energy density of the map part as
\begin{align}
e_k(\phi)=s^2|D^k\partial_t\phi|^2+|D^kD\phi|^2
\end{align}
and the $k$th energy as
\begin{align}
E_k(\phi)=\int_\Sigma e_k(\phi)\dv_{g_t}.
\end{align}
The $k$th total energy is given by
\begin{align}
F_k(\phi)=\sum_{l=0}^kE_l(\phi)=\int_\Sigma\sum_{l=0}^k e_l(\phi)\dv_{g_t}=s^2\|\nabla_t\phi\|^2_{H^k}+\|D\phi\|_{H^k}^2,
\end{align}
where the Sobolev norms are taken with respect to the metric $g_t$. In \eqref{totalenergy_spinor} below, we also define the total energy $F_k(\psi)$ of the spinor part, which we already need in the following proposition:
\begin{Prop}\label{prop-evol-energy-map}
Let $T>0$, $r\in\N$, $r>(n-1)/2$, $k\in \left\{1,\ldots,r\right\}$ and \((\phi,\psi)\) be a solution of \eqref{wave-phi}, \eqref{wave-psi} such that
\begin{align*}
&\phi\in C^0([0,T),H^{r+1}(\Sigma,P))\cap  C^1([0,T),H^r(\Sigma,P)), \\
&\psi\in C^0([0,T),H^{r}(M,SM\otimes\phi^\ast TP))\cap C^1([0,T),H^{r-1}(M,SM\otimes\phi^\ast TP)).
\end{align*}
Then the $k$th energy of the map satisfies the following inequality
\begin{equation}\begin{split}
\label{evolution-energy-map}
\frac{d}{dt}E_k(\phi)
&\leq C_1(k)\big(\|\dot g\|_{C^{k-1}}+\|\partial_t\log N\|_{C^k}+s^{-1}\|D\log N\|_{C^k}+s^{-1}\|R^{\Sigma}\|_{C^{k-1}}\big)F_k(\phi)\\&\quad-(n-2)s\dot s\int_{\Sigma}|D^k\nabla_t\phi|^2\dv_{g_t} \\
&\quad+ C_2(k,n,g)s^{-1}\|R^P\|_{C^{k-1}}\sum_{l=1}^{k-1}F_r(\phi)^{2+l/2}\\ 
&\quad+C_3(k,n,g)s^{2-n}\|N^{2-n}\|_{C^k}\|R^P\|_{C^k}\|\partial_t\|_{C^k}\sum_{l=0}^k F_r(\phi)^{1+l/2}F_{r-1}(\psi)\\
&\quad +C_4(k,n,g)s^{1-n}\|N^{2-n}\|_{C^k}\|R^P\|_{C^k}\sum_{l=0}^k  F_r(\phi)^{1+l/2}F_{r-1}(\psi)\\
&\quad+ C_5(k,n,g)s^{3-2n}\|N^{4-2n}\|_{C^k}\|R^P\|_{C^{k+1}}
\sum_{l=0}^k F_r(\phi)^{1/2+l/2}F_{r-1}(\psi)^2,
\end{split}
\end{equation}
where the positive constants $C_i$, $i\in\left\{1,\ldots,5\right\}$ depend on $n,k$ and the Sobolev constant of the metric $g_t$.
\end{Prop}
\begin{proof}
Assume for the moment that the initial data is compactly supported such that the solution is spacelike compactly supported.
Using the general formula \eqref{energy-evolution-general} we find
\begin{align*}
\frac{d}{dt}E_k(\phi)=&2\int_\Sigma\langle D^k\Box_h\phi,D^k\nabla_t\phi\rangle \dv_{g_t}
+\frac{1}{2}\int_\Sigma(|D^kD\phi|-s^2|D^k\nabla_t\phi|^2)\tr_g\dot g \dv_{g_t} \\
\nonumber &+2s^2\int_\Sigma\langle[\nabla_t,D^k]\nabla_t\phi,D^k\nabla_t\phi\rangle \dv_{g_t}
+2\int_\Sigma\langle[D^\ast D,D^k]\phi, D^k\nabla_t\phi\rangle \dv_{g_t} \\
\nonumber &+2\int_\Sigma\langle[\nabla_t,D^{k+1}]\phi,D^{k+1}\phi\rangle \dv_{g_t}.
\end{align*}
Due to finite speed of propagation and an exhaustion procedure, this equality also holds generally for solutions in the above space. 
Note also that we will use Lemma \ref{sobolev_multiplication} frequently in the proof without mentioning it explicitly.
We have to estimate all terms on the right hand side and start by estimating the commutator terms
\begin{align*}
&s^2\int_\Sigma\langle D^k\nabla_t\phi,[\nabla_t,D^k]\nabla_t\phi\rangle \dv_{g_t} \\
&=s^2\sum_{\sum l_i+\sum {m_j}=k-1}\int_\Sigma\nabla^{l_1}R^P\star\underbrace{D^{m_1+1}\phi\star\ldots\star D^{m_{l_1}+1}\phi}_{l_1-\textrm{times}}
\star D^{l_2}\partial_t\phi\star D^{l_3+1}\phi\star D^{l_4}\partial_t\phi\star D^k\partial_t\phi \dv_{g_t}\\
&\quad+s^2\sum_{l=0}^{k-1}\int_\Sigma D^{k-l}\dot g\star D^{l+1}\partial_t\phi\star D^k\partial_t\phi \dv_{g_t} \\
&\leq C(k)s^2\|R^P\|_{C^{k-1}}\sum_{l=0}^{k-1}\|D\phi\|^{l+1}_{H^r}\|\partial_t\phi\|^3_{H^r}
+C(k)s^2\|\dot g\|_{C^{k-1}}\|\partial_t\phi\|^2_{H^k}\\
&\leq C(k)s^{-1}\|R^P\|_{C^{k-1}}\sum_{l=0}^{k-1}(F_r(\phi))^{2+l/2}
+C(k)\|\dot g\|_{C^{k-1}}F_k(\phi).
\end{align*}

The second commutator can be controlled as follows
\begin{align*}
&\int_\Sigma\langle D^{k+1}\phi,[\nabla_t,D^{k+1}\phi]\rangle \dv_{g_t} \\
&=\sum_{\sum l_i+\sum {m_j}=k-1}\int_\Sigma\nabla^{l_1} R^P\star\underbrace{D^{m_1+1}\phi\star\ldots\star D^{m_{l_1}+1}\phi}_{l_1-\textrm{times}}
\star D^{l_2}\partial_t\phi\star D^{l_3}D\phi\star D^{l_4}D\phi\star D^{k}D\phi \dv_{g_t} \\
&\quad+\sum_{l=0}^{k}\int_\Sigma D^{k-l}\dot g\star D^{l}D\phi\star D^kD\phi \dv_{g_t} \\
&\leq C(k)\|R^P\|_{C^{k-1}}\sum_{l=0}^{k-1}\|D\phi\|_{H^r}^{l+3}\|\partial_t\phi\|_{H^r}
+C(k)\|\dot g\|_{C^{k-1}}\|D\phi\|^2_{H^k}\\
&\leq C(k)s^{-1}\|R^P\|_{C^{k-1}}\sum_{l=0}^{k-1}(F_r(\phi))^{2+l/2}
+C(k)\|\dot g\|_{C^{k-1}}F_k(\phi).
\end{align*}

The third commutator can be estimated as follows
\begin{align*}
&\int_\Sigma\langle D^k\partial_t\phi,[D^\ast D,D^k]\phi\rangle \dv_{g_t} \\
&=\sum_{l=0}^{k-1}\int_\Sigma D^lR^{\Sigma}\star D^{k-l}\phi\star D^k\partial_t\phi \dv_{g_t} \\
&\quad+\sum_{\sum l_i+\sum {m_j}=k-1}\int_\Sigma\nabla^{l_1}R^P\star\underbrace{D^{m_1+1}\phi\star\ldots\star D^{m_{l_1}+1}\phi}_{l_1-\textrm{times}} 
\star D^{l_2}D\phi\star D^{l_3}D\phi\star D^{l_4}D\phi\star D^k\nabla_t\phi \dv_{g_t} \\
&\leq C(k)\|R^{\Sigma}\|_{C^{k-1}}\|D\phi\|_{H^{k-1}}\|\partial_t\phi\|_{H^k}
+C(k)\|R^P\|_{C^{k-1}}\sum_{l=0}^{k-1}\|D\phi\|^{l+3}_{H^r}\|\partial_t\phi\|_{H^r}\\
&\leq C(k)s^{-1}\|R^{\Sigma}\|_{C^{k-1}}F_k(\phi)+C(k)s^{-1}\|R^P\|_{C^{k-1}}\sum_{l=0}^{k-1}(F_r(\phi))^{2+l/2}.
\end{align*}

As a second step we estimate the terms that arise when inserting the equation \eqref{wave-phi} for \(\phi\) into
\eqref{energy-evolution-general}. In order to estimate the first term we calculate
\begin{align*}
(Ns)^{-1}\nabla_{\nabla (Ns)}\phi=-s\dot s\nabla_t\phi-s^2N^{-1}\partial_tN\nabla_t\phi +N^{-1}\langle DN, D\phi\rangle_g.
\end{align*}

This allows us to derive the following estimate
\begin{align*}
\int_{\Sigma}&\langle D^k\big((Ns)^{-1}\nabla_{\nabla (Ns)}\phi\big),D^k\nabla_t\phi\rangle \dv_{g_t} \\
=&-s\dot s\int_{\Sigma}|D^k\nabla_t\phi|^2\dv_{g_t}
+\sum_{l=0}^k\int_{\Sigma} D^{l}\big(\frac{DN}{N}\big)\star D^{k-l}D\phi\star D^k\nabla_t\phi \dv_{g_t} \\
&+s^2\sum_{l_i}^k\int_{\Sigma} D^{l_1}\big(\frac{\partial_tN}{N}\big)\star D^{l_2}\nabla_t\phi\star D^k\nabla_t\phi\dv_{g_t} \\
\leq & -s\dot s\int_{\Sigma}|D^k\nabla_t\phi|^2\dv_{g_t}+C(k)s^{-1}\|D\log N\|_{C^k}F_k(\phi) +C(k)\|\partial_t\log N\|_{C^k}F_k(\phi).
\end{align*}

In order to treat the second term on the right hand side of \eqref{wave-phi} 
we have to consider spatial and time derivatives separately
\begin{align*}
D^k\big(R^P(\psi&,i\partial_t\cdot\psi)d\phi(\partial_t)\big) \\
&=\sum_{\sum l_i+\sum {m_j}=k} \nabla^{l_1}R^P\star\underbrace{D^{m_1+1}\phi\star\ldots\star D^{m_{l_1}+1}\phi}_{l_1-\textrm{times}}\star
D^{l_2}\partial_t\star D^{l_3}\psi\star D^{l_4}\psi\star D^{l_5}\nabla_t\phi, \\
D^k\big(R^P(\psi&,i\partial_i\cdot\psi)d\phi(\partial_i)\big) \\
&=\sum_{\sum l_i+\sum {m_j}=k} \nabla^{l_1}R^P\star\underbrace{D^{m_1+1}\phi\star\ldots\star D^{m_{l_1}+1}\phi}_{l_1-\textrm{times}}\star 
D^{l_2}\psi\star D^{l_3}\psi\star D^{l_4}D\phi,
\end{align*}
such that we find
\begin{align*}
D^k&h^{\alpha\beta}\big(N^{2-n}R^P(\psi,i\partial_\alpha\cdot\psi)d\phi(\partial_\beta)\big)\\
=&\sum_{\sum l_i+\sum {m_j}=k} D^{l_1}N^{2-n}\star\big(s^2\nabla^{l_2}R^P\star\underbrace{D^{m_1+1}\phi\star\ldots\star D^{m_{l_2}+1}\phi}_{l_2-\textrm{times}}
\star\nonumber D^{l_3}\partial_t\star D^{l_4}\psi\star D^{l_5}\psi\star D^{l_6}\nabla_t\phi\\
&\nonumber+\nabla^{l_2}R^P\star\underbrace{D^{m_1+1}\phi\star\ldots\star D^{m_{l_2}+1}\phi}_{l_2-\textrm{times}} 
\star D^{l_3}\psi\star D^{l_4}\psi\star D^{l_5}D\phi
\big).
\end{align*}
This allows us to derive the following estimate:
\begin{align*}
\int_{\Sigma}\langle&  D^k\big(-\frac{1}{2}h^{\alpha\beta}(Ns)^{2-n}R^P(\psi,i\partial_\alpha\cdot\psi)d\phi(\partial_\beta)\big),D^k\nabla_t\phi\rangle \dv_{g_t}\\
&=-\frac{s^{2-n}}{2}\sum_{\sum l_i+\sum {m_j}=k}\int_{\Sigma}D^{l_1}N^{2-n} 
\\ &\quad \star s^2\big(\nabla^{l_2}R^P\star\underbrace{D^{m_1+1}\phi\star\ldots\star D^{m_{l_2}+1}\phi}_{l_2-\textrm{times}}
\star\nonumber D^{l_3}\partial_t\star D^{l_4}\psi\star D^{l_5}\psi\star D^{l_6}\nabla_t\phi\star D^k\nabla_t\phi \\
\nonumber&\quad+D^{l_2}R^P\star\underbrace{D^{m_1+1}\phi\star\ldots\star D^{m_{l_2}+1}\phi}_{l_2-\textrm{times}} 
\star D^{l_3}\psi\star D^{l_4}\psi\star D^{l_5}D\phi\star D^k\nabla_t\phi
\big)\dv_{g_t}\\
&\leq C(k,n)s^{4-n}\|N^{2-n}\|_{C^k}\|R^P\|_{C^k}\|\partial_t\|_{C^k}\sum_{l=0}^k\|D\phi\|^l_{H^r}\|\psi\|^2_{H^r}\|\nabla_t\phi\|^2_{H^r} \\
&\quad+C(k,n)s^{2-n}\|N^{2-n}\|_{C^k}\|R^P\|_{C^k}\sum_{l=0}^k\|D\phi\|^{l+1}_{H^r}\|\psi\|^2_{H^r}\|\nabla_t\phi\|_{H^r}\\
&\leq C(k,n)s^{2-n}\|N^{2-n}\|_{C^k}\|R^P\|_{C^k}\|\partial_t\|_{C^k}\sum_{l=0}^k F_r(\phi)^{1+l/2}F_{r-1}(\psi)\\
&\quad+C(k,n)s^{1-n}\|N^{2-n}\|_{C^k}\|R^P\|_{C^k}\sum_{l=0}^k  F_r(\phi)^{1+l/2}F_{r-1}(\psi)
.
\end{align*}
The third term on the right hand side of \eqref{wave-phi} can be computed as
\begin{align*}
&D^k\big(N^{4-2n}\langle(\nabla R^P)^\sharp(\psi,\psi)\psi,\psi\rangle\big)\\
&=\sum_{\sum l_i+\sum {m_j}=k} 
D^{l_1}N^{4-2n}\star \nabla^{l_2+1}R^P\star\underbrace{D^{m_1+1}\phi\star\ldots\star D^{m_{l_2}+1}\phi}_{l_2-\textrm{times}}\star D^{l_3}\psi\star D^{l_4}\psi\star D^{l_5}\psi\star D^{l_6}\psi.
\end{align*}
This allows us to derive the following estimate
\begin{align*}
\int_{\Sigma}\langle D^k(&(Ns)^{4-2n}\langle(\nabla R^P)^\sharp(\psi,\psi)\psi,\psi\rangle),D^k\partial_t\phi\rangle \dv_{g_t}\\
=&s^{4-2n}\sum_{\sum l_i+\sum {m_j}=k}\int_\Sigma D^{l_1}N^{4-2n}\star \nabla^{l_2+1}R^P\\&\quad\star\underbrace{D^{m_1+1}\phi\star\ldots\star D^{m_{l_2}+1}\phi}_{l_2-\textrm{times}}\star D^{l_3}\psi\star D^{l_4}\psi\star D^{l_5}\psi\star D^{l_6}\psi
\star D^{k}\partial_t\phi \dv_{g_t} \\
\leq &s^{4-2n}C(k)\|N^{4-2n}\|_{C^k}\|R^P\|_{C^{k+1}}\sum_{l=0}^k\|D\phi\|^l_{H^r}\|\psi\|^4_{H^r}\|\partial_t\phi\|_{H^r}\\
\leq & s^{3-2n}C(k)\|N^{4-2n}\|_{C^k}\|R^P\|_{C^{k+1}}
\sum_{l=0}^k F_r(\phi)^{1/2+l/2}F_{r-1}(\psi)^2.
\end{align*}
Adding up the different contributions concludes the proof.
\end{proof}

\subsection{Energy of the spinor}
We define the $k$th energy density of the spinor part as
\begin{align}
e_k(\psi)=s^2|D^k\nabla_t\psi|^2+|D^kD\psi|^2
\end{align}
and the $k$th energy as
\begin{align}
E_k(\psi)=\int_\Sigma e_k(\psi)\dv_{g_t}.
\end{align}
The $k$th total energy is given by
\begin{align}\label{totalenergy_spinor}
F_k(\psi)=\sum_{l=0}^kE_l(\psi)+\|\psi\|_{L^2}^2=\int_\Sigma\sum_{l=0}^k e_l(\psi)\dv_{g_t}+\|\psi\|_{L^2}^2=s^2\|\nabla_t\psi\|_{H^k}^2+\|\psi\|_{H^{k+1}}^2.
\end{align}
Note that 
\begin{align}\label{L^2-spinor-inequality}
\frac{d}{dt}\|\psi\|_{L^2}^2\leq 2\|\nabla_t\psi\|_{L^2}\|\psi\|_{L^2}+\frac{1}{2}\|\trace\dot{g}\|_{L^{\infty}}\|\psi\|_{L^2}^2\leq (s^{-1}+\frac{1}{2}\|\trace\dot{g}\|_{L^{\infty}})F_1(\psi).
\end{align}

\begin{Prop}\label{prop-evol-energy-spinor}
Let $T>0$, $r\in\N$, $r>(n-1)/2$, $k\in \left\{1,\ldots,r-1\right\}$ and \((\phi,\psi)\) be a solution of \eqref{wave-phi}, \eqref{wave-psi} such that
\begin{align*}
&\phi\in C^0([0,T),H^{r+1}(\Sigma,P))\cap  C^1([0,T),H^r(\Sigma,P)), \\
&\psi\in C^0([0,T),H^{r}(M,SM\otimes\phi^\ast TP))\cap C^1([0,T),H^{r-1}(M,SM\otimes\phi^\ast TP)).
\end{align*}
Then the \(k\)th energy of the spinor satisfies the following inequality
\begin{equation}\begin{split}
\label{evolution-energy-spinor}
\frac{d}{dt}E_k(\psi)&\leq
C_1(k,n)(\|R^{SM}(\partial_t,\cdot)\|_{C^{k-1}}+\|\dot g\|_{C^k}+\|D^\ast e_0\|_{L^\infty})F_k(\psi)\\
&\quad+ C_2(n,k)\cdot s^{-1}(\|R^{\Sigma}\|_{C^{k-1}}+\|R^{SM}\|_{C^k}+\|\scal^M\|_{C^k})F_k(\psi)\\
&\quad+ C_3(n,k,g)(s^{-1}+\|\partial_t\|_{C^k})\|R^P\|_{C^{k+1}}\sum_{l=0}^kF_r(\phi)^{l/2+1}F_{r-1}(\psi)\\
&\quad+ C_4(n,k,g)s^{2-n}\|N^{2-n}\|_{C^k}\|R^P\|_{C^{k+1}}\|\partial_t\|_{C^k}\sum_{l=0}^{k+1}F_r(\phi)^{l/2}F_{r-1}(\psi)^2\\
&\quad+ (n-2)\cdot C_5(n,k,g)s^{3-n}\|N^{2-n}\partial_t\log N\|_{C^k}\|R^P\|_{C^{k}}\|\partial_t\|_{C^k}\sum_{l=0}^{k}F_r(\phi)^{l/2}F_{r-1}(\psi)^2\\
&\quad+ (n-2)\cdot C_6(n,k,g)\dot{s}s^{2-n}\|N^{2-n}\|_{C^k}\|R^P\|_{C^{k}}\|\partial_t\|_{C^k}\sum_{l=0}^{k}F_r(\phi)^{l/2}F_{r-1}(\psi)^2\\
&\quad+ C_7(n,k,g)s^{1-n}\|N^{2-n}\|_{C^{k+1}}\|R^P\|_{C^{k+1}}\sum_{l=0}^{k+1}F_r(\phi)^{l/2}F_{r-1}(\psi)^2\\
&\quad+ C_8(n,k,g)s^{1-n}\|N^{2-n}\|_{C^{k}}\|R^P\|_{C^{k}}^2\sum_{l=0}^{2k}F_r(\phi)^{l/2}F_{r-1}(\psi)^3,
\end{split}
\end{equation}
where the positive constants $C_i$, $i\in\left\{3,\ldots,8\right\}$ depend on $n,k$ and the Sobolev constant of the metric $g_t$.
\end{Prop}
\begin{Bem}
The curvature quantities coming from the spinor bundle above have to be understood as the following sections: 
$R^{SM}(\partial_t,\cdot)\in \Gamma(T^*\Sigma\otimes \mathrm{End}(SM))$, $R^{SM}=R^{SM}(.,.)\in \Gamma(\Lambda^2T\Sigma\otimes \mathrm{End}(SM))$, both by restricting to vectors tangent to $\Sigma$. 
Moreover $\partial_t\in\Gamma(\mathrm{End}(SM))$ acts by Clifford multiplication.
The $C^k$ norms with respect to $g_t$ are then defined in a canonical way.
\end{Bem}
\begin{proof}
In order to derive an energy estimate for the spinor we make use of \eqref{kenergy-spinor}.
As before, we may at first assume that the solution is spacelike compactly supported before we obtain this equality for the general case by an exhaustion procedure. 
We will also frequently use Lemma \ref{sobolev_multiplication} in this proof.
Again, we have to estimate all terms on the right hand side and 
start by estimating the commutator terms:
\begin{align*}
s^2\int_\Sigma&\langle D^k\nabla_t\psi,[\nabla_t,D^k]\psi\rangle \dv_{g_t} \\
&=s^2\sum_{l=0}^{k-1}\int_\Sigma D^l(R^{SM}(\partial_t,\cdot))\star D^{k-l-1}\nabla_t\psi\star D^k\nabla_t\psi \dv_{g_t} \\
&\quad+s^2\sum_{l=0}^{k-1}\int_\Sigma D^{k-l}\dot g\star D^{l+1}\nabla_t\psi\star D^k\nabla_t\psi \dv_{g_t}\\
&\quad+s^2\sum_{\sum l_i+\sum {m_j}=k-1}\int_\Sigma\nabla^{l_1} R^P\\
&\qquad\star\underbrace{D^{m_1+1}\phi\star\ldots\star D^{m_{l_1}+1}\phi}_{l_1-\textrm{times}}
\star D^{l_2}\nabla_t\phi\star D^{l_3}D\phi\star D^{l_4}\nabla_t\psi\star D^k\nabla_t\psi \dv_{g_t}\\
&\leq s^2 C(k)\|R^{SM}(\partial_t,\cdot)\|_{C^{k-1}}\|\nabla_t\psi\|^2_{H^k}
+s^2C(k)\|\dot g\|_{C^k}\|\nabla_t\psi\|^2_{H^k} \\
&\qquad+s^2C(k)\|R^P\|_{C^{k-1}}\sum_{l=0}^{k-1}\|D\phi\|_{H^r}^{l+1}\|\partial_t\phi\|_{H^r}\|\nabla_t\psi\|^2_{H^{r-1}}\\
&\leq C(k)\|R^{SM}(\partial_t,\cdot)\|_{C^{k-1}} F_k(\psi)+
C(k)\|\dot g\|_{C^k} F_k(\psi)\\
&\qquad+s^{-1}C(k)\|R^P\|_{C^{k-1}}\sum_{l=0}^{k-1}F_r(\phi)^{l/2+1}F_{r-1}(\psi).
\end{align*}
The second commutator term can be estimated as follows:
\begin{align*}
&\int_\Sigma\langle[D^\ast D,D^k]\psi, D^k\nabla_t\psi\rangle \dv_{g_t} \\
&=\sum_{l=0}^{k-1}\int_\Sigma D^lR^M\star D^{k-l}\psi\star D^k\nabla_t\psi \dv_{g_t}
+\sum_{l=0}^k\int_\Sigma D^l R^{SM}\star D^{k-l}\psi\star D^k\nabla_t\psi \dv_{g_t} \\
&\quad+\sum_{\sum l_i+\sum {m_j}=k}\int_\Sigma\nabla^{l_1} R^P\star\underbrace{D^{m_1+1}\phi\star\ldots\star D^{m_{l_1}+1}\phi}_{l_1-\textrm{times}}
\star D^{l_2}D\phi\star D^{l_3}D\phi\star D^{l_4}\psi\star D^k\nabla_t\psi \dv_{g_t} \\
&\leq C(k)\|R^M\|_{C^{k-1}}\|\psi\|_{H^k}\|\nabla_t\psi\|_{H^k}
+C(k)\|R^{SM}\|_{C^{k}}\|\psi\|_{H^{k}}\|\nabla_t\psi\|_{H^k} \\
&\quad+C(k)\|R^P\|_{C^{k}}\sum_{l=0}^{k}\|D\phi\|_{H^r}^{l+2}\|\psi\|_{H^r}\|\nabla_t\psi\|_{H^{r-1}}\\
&\leq s^{-1}C(k)(\|R^M\|_{C^{k-1}}+\|R^{SM}\|_{C^{k}})F_k(\psi)+s^{-1}C(k)\|R^P\|_{C^{k}}\sum_{l=0}^{k}F_r(\phi)^{l/2+1}F_{r-1}(\psi).
\end{align*}
The third commutator can be controlled as follows:
\begin{align*}
&\int_\Sigma\langle[\nabla_t,D^{k+1}]\psi,D^{k+1}\psi\rangle \dv_{g_t}\\
&=\sum_{l=0}^k\int_\Sigma D^l(R^{SM}(\partial_t,\cdot))\star D^{k-l}\psi\star D^{k+1}\psi \dv_{g_t}
+\sum_{l=0}^k\int_\Sigma D^{k-l-1}\dot g\star D^{l+1}\psi\star D^{k+1}\psi \dv_{g_t} \\
&\quad+\sum_{\sum l_i+\sum {m_j}=k-1}\int_\Sigma\nabla^{l_1} R^P\star\underbrace{D^{m_1+1}\phi\star\ldots\star D^{m_{l_1}+1}\phi}_{l_1-\textrm{times}}
\star D^{l_2}\partial_t\phi\star D^{l_3}D\phi\star D^{l_4}\psi\star D^{k+1}\psi \dv_{g_t} \\
&\leq  C(k)\|R^{SM}(\partial_t,\cdot)\|_{C^k}\|\psi\|^2_{H^{k+1}}+C(k)\|\dot g\|_{C^k}\|\psi\|^2_{H^{k+1}} \\
&\quad+C(k)\|R^P\|_{C^k}\sum_{l=0}^{k}\|D\phi\|_{H^r}^{l+1}\|\nabla_t\phi\|_{H^r}\|\psi\|^2_{H^{r}}\\
&\leq C(k)\|R^{SM}(\partial_t,\cdot)\|_{C^k}F_k(\psi)+
C(k)\|\dot g\|_{C^k}F_k(\psi)
+s^{-1}C(k)\|R^P\|_{C^k}\sum_{l=0}^{k}F_r(\phi)^{l/2+1}F_{r-1}(\psi).
\end{align*}
As a next step we estimate the terms that arise when inserting \eqref{wave-psi}
into \eqref{energy-evolution-general}.
The first term can easily be controlled as
\begin{align*}
\int_\Sigma\langle D^k(\frac{\scal^M}{4}\psi),D^k\nabla_t\psi\rangle \dv_{g_t}=&\frac{1}{4}\sum_{l=0}^k\int_\Sigma D^{l}\scal^M\star D^{k-l}\psi\star D^k\nabla_t\psi \dv_{g_t} \\
\leq& C(k)\|\scal^M\|_{C^k}s^{-1}F_k(\psi).
\end{align*}
Regarding the second term on the right hand side of \eqref{wave-psi} 
we again expand space and time contributions
\begin{align*}
h^{\alpha\gamma}h^{\beta\delta}\partial_\alpha\cdot\partial_\beta\cdot R^P(d\phi(\partial_\gamma),d\phi(\partial_\delta))\psi
=&-2s^2g^{ij}\partial_t\cdot\partial_i\cdot R^P(d\phi(\partial_t),d\phi(\partial_j))\psi \\
&+g^{ij}g^{kl}\partial_i\cdot\partial_k\cdot R^P(d\phi(\partial_j),d\phi(\partial_l))\psi.
\end{align*}
Note that we do not get a term proportional to \(|d\phi(\partial_t)|^2\) due to
symmetry reasons.
We then find
\begin{align*}
D^k\big(&h^{\alpha\gamma}h^{\beta\delta}\partial_\alpha\cdot\partial_\beta\cdot R^P(d\phi(\partial_\gamma),d\phi(\partial_\delta))\psi)\\
=&\sum_{\sum l_i+\sum {m_j}=k}\big(D^{l_1}R^P\star\underbrace{D^{m_1+1}\phi\star\ldots\star D^{m_{l_1}+1}\phi}_{l_1-\textrm{times}} 
\star D^{l_2}D\phi\star D^{l_3}D\phi\star D^{l_4}\psi\big) \\
&+s^2\sum_{\sum l_i+\sum {m_j}=k}\big(D^{l_1}\partial_t\star D^{l_2}R^P\star\underbrace{D^{m_1+1}\phi\star\ldots\star D^{m_{l_2}+1}\phi}_{l_2-\textrm{times}} 
\star D^{l_3}d\phi(\partial_t)\star D^{l_4}d\phi\star D^{l_5}\psi\big).
\end{align*}
This allows us to derive the following estimate
\begin{align*}
\int_\Sigma&\langle D^k\big(h^{\alpha\gamma}h^{\beta\delta}\partial_\alpha\cdot\partial_\beta\cdot R^P(d\phi(\partial_\gamma),d\phi(\partial_\delta))\psi\big),D^k\nabla_t\psi\rangle \dv_{g_t}\\
&=\sum_{\sum l_i+\sum {m_j}=k}\int_\Sigma D^{l_1}R^P\star\underbrace{D^{m_1+1}\phi\star\ldots\star D^{m_{l_1}+1}\phi}_{l_1-\textrm{times}} 
\star D^{l_2}d\phi\star D^{l_3}d\phi\star D^{l_4}\psi\star D^k\nabla_t\psi \dv_{g_t} \\
&\quad+s^2\sum_{\sum l_i+\sum {m_j}=k}\int_\Sigma D^{l_1}\partial_t\star D^{l_2}R^P\star\underbrace{D^{m_1+1}\phi\star\ldots\star D^{m_{l_2}+1}\phi}_{l_2-\textrm{times}} \\& \qquad
\star D^{l_3}d\phi(\partial_t)\star D^{l_4}d\phi\star D^{l_5}\psi\star D^k\nabla_t\psi \dv_{g_t} \\
&\leq C(k)\|R^P\|_{C^k}\sum_{l=0}^k\|D\phi\|_{H^r}^{l+2}\|\psi\|_{H^{r}}\|\nabla_t\psi\|_{H^{k}} \\
&\quad+s^2C(k)\|\partial_t\|_{C^k}\|R^P\|_{C^k}\sum_{l=0}^k\|D\phi\|_{H^r}^{l+1}\|\nabla_t\phi\|_{H^r}\|\psi\|_{H^r}\|\nabla_t\psi\|_{H^{k}}\\
&\leq s^{-1}C(k)\|R^P\|_{C^k}\sum_{l=0}^k F_r(\phi)^{l/2+1}F_{r-1}(\psi)+C(k)\|\partial_t\|_{C^k}\|R^P\|_{C^k}\sum_{l=0}^k F_r(\phi)^{l/2+1}F_{r-1}(\psi).
\end{align*}
To manipulate the third term on the right hand side of \eqref{wave-psi} we first consider the derivative with respect to \(t\).
We find
\begin{align*}
\nabla_{t}(R^P(\psi,\psi)(Ns)^{2-n})\partial_t\cdot\psi=&(\nabla_{d\phi(\partial_t)}R^P)(\psi,\psi)(Ns)^{2-n}\partial_t\cdot\psi 
+2R^P(\nabla_t\psi,\psi)(Ns)^{2-n}\partial_t\cdot\psi\\
&+(2-n)R^P(\psi,\psi)(Ns)^{2-n}\partial_t\log N\partial_t\cdot\psi \\
&+(2-n)\dot ss^{1-n}N^{2-n}R^P(\psi,\psi)\partial_t\cdot\psi.
\end{align*}
As a next step we calculate the \(k\)-th spatial derivative of this expression
\begin{align*}
D^k\big(&(\nabla_{d\phi(\partial_t)}R^P)(\psi,\psi)N^{2-n}\partial_t\cdot\psi\big)=\sum_{\sum l_i+\sum {m_j}=k}^kD^{l_1}N^{2-n}\star D^{l_2+1}R^P\\
&\quad \star\underbrace{D^{m_1+1}\phi\star\ldots\star D^{m_{l_2+1}+1}\phi}_{l_2+1-\textrm{times}}
\star D^{l_3}\psi\star D^{l_4}\psi\star D^{l_5}\partial_t\star D^{l_6}\psi\star D^{l_7}d\phi(\partial_t),\\ 
D^k\big(&R^P(\nabla_t\psi,\psi)N^{2-n}\partial_t\cdot\psi\big)=\sum_{\sum l_i+\sum {m_j}=k}^kD^{l_1}N^{2-n}\star D^{l_2}R^P
\\&\star \underbrace{D^{m_1+1}\phi\star\ldots\star D^{m_{l_2}+1}\phi}_{l_2-\textrm{times}}
\star D^{l_3}\nabla_t\psi\star D^{l_4}\psi\star D^{l_5}\partial_t\star D^{l_6}\psi,\\
D^k\big(&R^P(\psi,\psi)N^{2-n}\partial_t\log N\partial_t\cdot\psi)=\sum_{\sum l_i+\sum {m_j}=k}^kD^{l_1}N^{2-n}\partial_t\log N\star D^{l_2}R^P\\
&\star \underbrace{D^{m_1+1}\phi\star\ldots\star D^{m_{l_2}+1}\phi}_{l_2-\textrm{times}}
\star D^{l_3}\psi\star D^{l_4}\psi\star D^{l_5}\partial_t\star D^{l_6}\psi,\\
D^k\big(&N^{2-n}R^P(\psi,\psi)\partial_t\cdot\psi\big)=\sum_{\sum l_i+\sum {m_j}=k}^kD^{l_1}N^{2-n}\star D^{l_2}R^P\\
&\star \underbrace{D^{m_1+1}\phi\star\ldots\star D^{m_{l_2}+1}\phi}_{l_2-\textrm{times}}
\star D^{l_3}\psi\star D^{l_4}\psi\star D^{l_5}\partial_t\star D^{l_6}\psi.
\end{align*}
These manipulations allow us to derive the following estimates:
\begin{align*}
s^2\int_\Sigma&\langle D^k\big(\nabla_{t}(R^P(\psi,\psi)(Ns)^{2-n})\partial_t\cdot\psi\big),D^k\nabla_t\psi\rangle \dv_{g_t} \\
=&s^{4-n}\sum_{\sum l_i+\sum {m_j}=k}\int_\Sigma D^{l_1}N^{2-n}\star D^{l_2+1}R^P\star \underbrace{D^{m_1+1}\phi\star\ldots\star D^{m_{l_2+1}+1}\phi}_{l_2+1-\textrm{times}}\\&\quad
\star D^{l_3}\psi\star D^{l_4}\psi\star D^{l_5}\partial_t\star D^{l_6}\psi\star D^{l_7}d\phi(\partial_t)\star D^k\nabla_t\psi \dv_{g_t} \\
&+s^{4-n}\sum_{\sum l_i+\sum {m_j}=k}\int_\Sigma D^{l_1}N^{2-n}\star D^{l_2}R^P\star \underbrace{D^{m_1+1}\phi\star\ldots\star D^{m_{l_2}+1}\phi}_{l_2-\textrm{times}}\\&\quad
\star D^{l_3}\nabla_t\psi\star D^{l_4}\psi\star D^{l_5}\partial_t\star D^{l_6}\psi\star D^k\nabla_t\psi \dv_{g_t} \\
&+(2-n)\dot ss^{3-n}\sum_{\sum l_i+\sum {m_j}=k}\int_\Sigma D^{l_1}N^{2-n}\star D^{l_2}R^P\star \underbrace{D^{m_1+1}\phi\star\ldots\star D^{m_{l_2}+1}\phi}_{l_2-\textrm{times}}\\&\quad
\star D^{l_3}\psi\star D^{l_4}\psi\star D^{l_5}\partial_t\star D^{l_6}\psi\star D^k\nabla_t\psi \dv_{g_t} \\
&+(2-n)s^{4-n}\sum_{\sum l_i+\sum {m_j}=k}\int_\Sigma D^{l_1}N^{2-n}\partial_t\log N\star D^{l_2}R^P\star \underbrace{D^{m_1+1}\phi\star\ldots\star D^{m_{l_2}+1}\phi}_{l_2-\textrm{times}}\\&\quad
\star D^{l_3}\psi\star D^{l_4}\psi\star D^{l_5}\partial_t\star D^{l_6}\psi\star D^k\nabla_t\psi \dv_{g_t} \\
\leq& C(k)s^{4-n}\|N^{2-n}\|_{C^k}\|R^P\|_{C^{k+1}}\|\partial_t\|_{C^k}\sum_{l=0}^k\|D\phi\|^l_{H^r}\|\psi\|^3_{H^r}\|\nabla_t\phi\|_{H^r}\|\nabla_t\psi\|_{H^k} \\
&+C(k)s^{4-n}\|N^{2-n}\|_{C^k}\|R^P\|_{C^{k}}\|\partial_t\|_{C^k}\sum_{l=0}^k\|D\phi\|^l_{H^r}\|\psi\|^2_{H^r}\|\nabla_t\psi\|^2_{H^k} \\
&+C(k,n)(n-2)\dot ss^{3-n}\|N^{2-n}\|_{C^k}\|R^P\|_{C^{k}}\|\partial_t\|_{C^k}\sum_{l=0}^k\|D\phi\|^l_{H^r}\|\psi\|^3_{H^r}\|\nabla_t\psi\|_{H^k}\\
&+C(k,n)(n-2)s^{4-n}\|N^{2-n}\partial_t\log N\|_{C^k}\|R^P\|_{C^{k}}\|\partial_t\|_{C^k}\sum_{l=0}^k\|D\phi\|^l_{H^r}\|\psi\|^3_{H^r}\|\nabla_t\psi\|_{H^k}\\
\leq& C(k)s^{2-n}\|N^{2-n}\|_{C^k}\|R^P\|_{C^{k+1}}\|\partial_t\|_{C^k}\sum_{l=0}^kF_r(\phi)^{l/2+1/2}F_{r-1}(\psi)^2\\
&+C(k)s^{2-n}\|N^{2-n}\|_{C^k}\|R^P\|_{C^{k}}\|\partial_t\|_{C^k}
\sum_{l=0}^kF_r(\phi)^{l/2}F_{r-1}(\psi)^2\\ 
&+C(k,n)(n-2)\dot ss^{2-n}\|N^{2-n}\|_{C^k}\|R^P\|_{C^{k}}\|\partial_t\|_{C^k}\sum_{l=0}^kF_r(\phi)^{l/2}F_{r-1}(\psi)^2
\\
&+C(k,n)(n-2)s^{3-n}\|N^{2-n}\partial_t\log N\|_{C^k}\|R^P\|_{C^{k}}\|\partial_t\|_{C^k}\sum_{l=0}^kF_r(\phi)^{l/2}F_{r-1}(\psi)^2.
\end{align*}
As a next step we take care of the spatial derivatives in the third term on the right hand side of
\eqref{wave-psi}. These can be computed as
\begin{align*}
D(R^P(\psi,\psi)(Ns)^{2-n})\cdot\psi=&(D_{d\phi}R^P)(\psi,\psi)(Ns)^{2-n}\cdot\psi
+2R^P(D\psi,\psi)(Ns)^{2-n}\cdot\psi \\
&+s^{2-n}R^P(\psi,\psi)D(N)^{2-n}\cdot\psi.
\end{align*}
The \(k\)th spatial derivative of these terms acquire the forms
\begin{align*}
D^k\big((&D_{d\phi}R^P)(\psi,\psi)(Ns)^{2-n}\cdot\psi\big)\\
&=s^{2-n}\sum_{\sum l_i+\sum {m_j}=k} D^{l_1+1}R^P\star\underbrace{D^{m_1+1}\phi\star\ldots\star D^{m_{l_1+1}+1}\phi}_{(l_1+1)-\textrm{times}}
\star D^{l_2}\psi\star D^{l_3}\psi\star D^{l_4}\psi\star D^{l_5}N^{2-n}, \\
D^k\big(&R^P(D\psi,\psi)(Ns)^{2-n}\cdot\psi\big)\\
&=s^{2-n}\sum_{\sum l_i+\sum {m_j}=k} D^{l_1}R^P\star\underbrace{D^{m_1+1}\phi\star\ldots\star D^{m_{l_1}+1}\phi}_{l_1-\textrm{times}}
\star D^{l_2}D\psi\star D^{l_3}\psi\star D^{l_4}\psi\star D^{l_5}N^{2-n}, \\
D^k\big(&(s)^{2-n}R^P(\psi,\psi)D(N)^{2-n}\cdot\psi\big)\\
&=s^{2-n}\sum_{\sum l_i+\sum {m_j}=k} D^{l_1}R^P\star\underbrace{D^{m_1+1}\phi\star\ldots\star D^{m_{l_1}+1}\phi}_{l_1-\textrm{times}}
\star D^{l_2}\psi\star D^{l_3}\psi\star D^{l_4}\psi\star D^{l_5}DN^{2-n}.
\end{align*}
These manipulations allow us to derive the following estimates:
\begin{align*}
\int_\Sigma\langle& D^k\big(D(R^P(\psi,\psi)(Ns)^{2-n})\partial_i\cdot\psi\big),D^k\nabla_t\psi\rangle \dv_{g_t} \\
&=s^{2-n}\sum_{\sum l_i+\sum {m_j}=k}\int_\Sigma D^{l_1}N^{2-n}\star D^{l_2+1}R^P\star \underbrace{D^{m_1+1}\phi\star\ldots\star D^{m_{l_2+1}+1}\phi}_{(l_2+1)-\textrm{times}}\\&\quad
\star D^{l_3}\psi\star D^{l_4}\psi\star D^{l_5}\psi\star D^k\nabla_t\psi \dv_{g_t} \\
&\quad+s^{2-n}\sum_{\sum l_i+\sum {m_j}=k}\int_\Sigma D^{l_1}N^{2-n}\star D^{l_2}R^P\star \underbrace{D^{m_1+1}\phi\star\ldots\star D^{m_{l_2}+1}\phi}_{l_2-\textrm{times}}\\&\quad
\star D^{l_3}D\psi\star D^{l_4}\psi\star D^{l_5}\psi\star D^k\nabla_t\psi \dv_{g_t} \\
&\quad+s^{2-n}\sum_{\sum l_i+\sum {m_j}=k}\int_\Sigma D^{l_1}D(N^{2-n})\star D^{l_2}R^P\star \underbrace{D^{m_1+1}\phi\star\ldots\star D^{m_{l_2}+1}\phi}_{l_2-\textrm{times}}\\&\quad
\star D^{l_3}\psi\star D^{l_4}\psi\star D^{l_5}\psi\star D^k\nabla_t\psi \dv_{g_t} \\
&\leq C(k)s^{2-n}\|N^{2-n}\|_{C^k}\|R^P\|_{C^{k+1}}\sum_{l=0}^k\|D\phi\|^{l+1}_{H^r}\|\psi\|^3_{H^r}\|\nabla_t\psi\|_{H^k} \\
&\quad+C(k)s^{2-n}\|N^{2-n}\|_{C^k}\|R^P\|_{C^{k}}\sum_{l=0}^k\|D\phi\|^{l}_{H^r}\|\psi\|^3_{H^r}\|\nabla_t\psi\|_{H^k} \\
&\quad+C(k)s^{2-n}\|N^{2-n}\|_{C^{k+1}}\|R^P\|_{C^{k}}\sum_{l=0}^k\|D\phi\|^l_{H^r}\|\psi\|^3_{H^r}\|\nabla_t\psi\|_{H^k}\\
&\leq C(k)s^{1-n}\|N^{2-n}\|_{C^{k+1}}\|R^P\|_{C^{k+1}}\sum_{l=0}^{k+1}F_r(\phi)^{l/2}F_{r-1}(\psi)^2.
\end{align*}
To control the last term from the right hand side of \eqref{wave-psi} we calculate
\begin{align*}
D^k&\big((Ns)^{2-n}R^P(\psi,\psi)R^P(\psi,\psi)\psi\big) \\
=&s^{2-n}\sum_{\sum l_{i_1}+\sum m_{i_2}+\sum q_{i_3}=k} D^{l_1}N^{2-n}\star
D^{l_2}R^P\star\underbrace{D^{m_1+1}\phi\star\ldots\star D^{m_{l_2+1}+1}\phi}_{l_2-\textrm{times}} 
\\ &\star D^{l_3}R^P\star\underbrace{D^{q_1+1}\phi\star\ldots\star D^{q_{l_3+1}+1}\phi}_{l_3-\textrm{times}}\star D^{l_4}\psi\star D^{l_5}\psi\star D^{l_6}\psi\star D^{l_7}\psi\star D^{l_8}\psi.
\end{align*}
Consequently, we obtain the following estimate
\begin{align*}
\int_\Sigma&\langle D^k((Ns)^{2-n}R^P(\psi,\psi)R^P(\psi,\psi)\psi),D^k\nabla_t\psi\rangle \dv_{g_t} \\
&\leq s^{2-n}C(k)\|N^{2-n}\|_{C^k}\|R^P\|^2_{C^k}\sum_{l=0}^{2k}\|D\phi\|^{l}_{H^r}\|\psi\|^5_{H^r}\|\nabla_t\psi\|_{H^k}\\
&\leq s^{1-n}C(k)\|N^{2-n}\|_{C^k}\|R^P\|^2_{C^k}\sum_{l=0}^{2k} F_r(\phi)^{l/2}F_{r-1}(\psi)^3.
\end{align*}
The last term of \eqref{kenergy-spinor} can be estimated as
\begin{align*}
\int_{\Sigma}\langle D^{k}\psi,e_0\cdot(D^\ast e_0)\cdot D^{k+1}\psi\rangle\dv_{g_t}\leq
\|D^\ast e_0\|_{L^\infty}\|\psi\|^2_{H^r}.
\end{align*}

Adding up the different contributions yields the claim.
\end{proof}
\begin{Bem}
Note that all terms on the right hand side of \eqref{evolution-energy-spinor} have a similar analytic
structure except the terms proportional to \(C_3\) and \(C_8\), which contain higher powers of the \(H^r\)-norms
of \(\phi\) and \(\psi\).
\end{Bem}

\subsection{Energy of the coupled pair}
Before we prove an energy estimate for the coupled pair, we give more geometric interpretations of the terms appearing in the above estimates.
At first, recall that the second fundamental form $\sff\in \Gamma(T^*\Sigma^{\odot2})$ of a hypersurface $(\left\{t\right\}\times\Sigma, g_t)\subset (M,h)$ is given by
\begin{align*}
\sff(X,Y)=\langle \nabla_XY-D_XY,\nu\rangle_h=\frac{1}{2}\langle s^2\dot{g}(X,Y)\cdot \partial_t,\nu\rangle_h=-\frac{s}{2}\dot{g}(X,Y),
\end{align*}
where $\nu$ is the future-directed unit normal of the hypersurface. 
\begin{Lem}\label{geometric_estimates}
We have the estimates
\begin{align*}
\|\dot{g}\|_{C^k}&\leq 2s^{-1} \|\sff\|_{C^k},\qquad \|De_0\|_{L^\infty}\leq\|\sff\|_{L^\infty},\qquad 
\|\partial_t\|_{C^k}\leq s^{-1}(1+\|\sff\|_{C^{k-1}}),\\
 \|R^{\Sigma}\|_{C^k}&\leq C(k)\sum_{l=0}^k\left\|\sff\right\|_{C^k}^l\|R^M\|_{C^k}+C(k)\left\|\sff\right\|_{C^k}^2, \\
 \|R^{SM}\|_{C^k}&\leq C(k,n)\sum_{l=0}^k\left\|\sff\right\|_{C^k}^l\|R^M\|_{C^k}, \\
\|R^{SM}(\partial_t,\cdot)\|_{C^k}&\leq C(k,n)s^{-1}\sum_{l=0}^{k+1}\left\|\sff\right\|^l_{C^k}\|R^M\|_{C^k}.
\end{align*}
Here, we defined the $C^k$-norm of $R^M$ by taking the $k$th covariant derivative with respect to $h$ but taking the norm with respect to the Riemannian reference metric $s^{-2}dt^2+g_t$ in order to get a nonnegative quantity.
\end{Lem}
\begin{proof}
The first estimate follows from the definition of $\sff$, the second and third from the facts that $\nabla_X\partial_t=\frac{1}{2}\dot{g}(X)$ and \(e_0=s\partial_t\). 
Now we recall that for every tensor $T\in \Gamma(T^*\Sigma^{\otimes k})$ on the manifold $M$, the difference between the covariant derivatives $\nabla$ and $D$ can be expressed as
\begin{align*}
DT=\nabla T+\sff\star T,
\end{align*}
which by induction yields
\begin{align*}
D^kT=\sum_{l_1+l_2+\sum m_i=k} \underbrace{D^{m_1}\sff\star\ldots\star D^{m_{l_1}}\sff}_{l_1-\text{times}}\star\nabla^{l_2}T.
\end{align*}
This formula in combination with the Gau{\ss} equation
\begin{align*}
R^{\Sigma}=R^{M}+\sff\star\sff
\end{align*}
yields the fourth inequality.
To prove the last two formulas, we recall that
\begin{align*}
R^{SM}(X,Y)\psi=\frac{1}{4}\sum_{\alpha,\beta=0}^{n-1}R^M(X,Y,\partial_{\alpha},\partial_{\beta})\partial_{\alpha}\cdot \partial_{\beta}\cdot \psi
\end{align*}
for $X,Y\in \Gamma(TM)$,
where $\left\{\partial_{\alpha}\right\}$ is a local pseudo-orthonormal frame.
Therefore, we get in the above situation
\begin{align*}
(D^k R^{SM})\psi=\sum_{l_1+l_2+\sum m_i=k} \underbrace{D^{m_1}\sff\star\ldots\star D^{m_{l_1}}\sff}_{l_1-\text{times}}\star\nabla^{l_2}R^M\star\psi,
\end{align*}
which yields the fifth inequality.
Similarly, by using the product rule, we obtain
\begin{align*}
(D^k R^{SM}(\partial_t,\cdot))\psi=\sum_{\sum l_i+\sum m_j=k} D^{l_1}\partial_t\star\underbrace{D^{m_1}\sff\star\ldots\star D^{m_{l_2}}\sff}_{l_2-\text{times}}\star\nabla^{l_3}R^M\star\psi
\end{align*}
and by using the third inequality we obtain the last one.
\end{proof}

At this point we are ready to control the total energy of \((\phi,\psi)\), which we define by
\begin{align*}
F_r(\phi,\psi)&=F_r(\phi)+F_{r-1}(\psi)=\sum_{k=0}^rE_k(\phi)+\sum_{k=0}^{r-1}E_k(\psi)+\|\psi\|_{L^2}^2 \\
&=s^2\|\partial_t\phi\|^2_{H^r}+\|D\phi\|^2_{H^r}+s^2\|\nabla_t\psi\|^2_{H^{r-1}}+\|D\psi\|^2_{H^{r-1}}+\|\psi\|_{L^2}^2.
\end{align*}

\begin{Prop}
Let $T>0$, $r\in\N$, $r>(n-1)/2$ and \((\phi,\psi)\) be a solution of \eqref{wave-phi}, \eqref{wave-psi} such that
\begin{align*}
&\phi\in C^0([0,T),H^{r+1}(\Sigma,P))\cap  C^1([0,T),H^r(\Sigma,P)), \\
&\psi\in C^0([0,T),H^{r}(M,SM\otimes\phi^\ast TP))\cap C^1([0,T),H^{r-1}(M,SM\otimes\phi^\ast TP)).
\end{align*}
Suppose that $\dot{s}\geq0$ and that the following uniform bounds
\begin{align*}
&0<C_2\leq N\leq C_3,\qquad \|N\|_{C^{r+1}}<C_4,\qquad \|\nabla_{\nu}N\|_{C^{r+1}}<C_5,\\
&\|\sff\|_{C^r}<C_6,\qquad \|R^M\|_{C^r}<C_7,\quad\|R^P\|_{C^{r+1}}<C_8, \quad \|\sff\|_{L^\infty}< C_9s^{-1}
\end{align*}
hold for some positive constants $C_i$, $i=1,\ldots,9$. Suppose finally that there is a uniform bound on all the Sobolev constants of $g_t$.
Then the total energy satisfies the following inequality
\begin{align}
\label{total-energy-inequality}
\frac{d}{dt}F_r(\phi,\psi)\leq  C\cdot s^{-1}\sum_{l=0}^{2r+4}F_r(\phi,\psi)^{l/2+1}
+C(n-2) \dot{s}s^{1-n}\sum_{l=0}^{r}F_r(\phi,\psi)^{l/2+2},
\end{align}
where $C>0$ depends on the above bounds.
\end{Prop}
\begin{proof}
The estimate is a direct consequence of the evolution inequalities \eqref{evolution-energy-map}, \eqref{L^2-spinor-inequality}, \eqref{evolution-energy-spinor}, Lemma \ref{geometric_estimates} and elementary estimates.
\end{proof}

\begin{proof}[Proof of Theorem \ref{result_diracwavemap}]
The existence of a local solution to the system \eqref{phi-dwmap-quer}, \eqref{psi-dwmap-quer} 
can be obtained by standard methods, see for example \cite[Proposition 9.12]{MR2527641}.
In order to prove long-time existence, it suffices to prove longtime existence of the solution $(\phi,\psi)$
of the equivalent system \eqref{wave-phi}, \eqref{wave-psi} on the conformal manifold $(M,h)$.

Moreover, by the continuation criterion for hyperbolic partial differential equations \cite[Lemma 9.14]{MR2527641}
if suffices to obtain a uniform bound on \(F_r(\phi,\psi)\) for all times \(t\in[0,\infty)\).

We set $f(t)=(n-2) s^{1-n}\dot{s}$. 
Note that \(f(t)\) is integrable with respect to \(t\). For $n=2$, this is trivial. For $n>2$, we get
due to the assumptions on $s$ that
\begin{align*}
(n-2)\int_0^\infty\dot s s^{1-n}dt=-\int_0^\infty \frac{d}{dt}(s^{2-n})dt
=-\big(s^{2-n}|_{t=\infty}-s^{2-n}|_{t=0}\big)=s(0)^{2-n}<\infty.
\end{align*}
As long as \(F_r(\phi,\psi)\leq 1\) we have the differential inequality
\begin{align*}
\frac{d}{dt}F_r(\phi,\psi)\leq C\big(s^{-1}(t)+f(t)\big)F_r(\phi,\psi),
\end{align*}
which can easily be integrated as
\begin{align*}
F_r(\phi,\psi)|_{t=T}\leq F_r(\phi,\psi)|_{t=0}\exp\big(\int_0^T(s^{-1}(t)+f(t))dt\big).
\end{align*}
Now, we set \(\Phi:=\int_0^\infty(s^{-1}(t)+f(t))dt<\infty\)
and choose \(\epsilon>0\) small enough such that \(F_r(\phi,\psi)|_{t=0}\leq (2\Phi)^{-1}\).
Suppose that \(T_0\) is the first time for which \(F_r(\phi,\psi)|_{T_0}=1\).
However, the energy inequality from above gives
\begin{align*}
F_r(\phi,\psi)|_{t=T_0}\leq \Phi F_r(\phi,\psi)|_{t=0}=\frac{1}{2}<1,
\end{align*}
which yields a contradiction.
Therefore we can conclude that
\begin{align*}
F_r(\phi,\psi)<1<\infty
\end{align*}
for all times \(t\in[0,\infty)\) completing the proof.
\end{proof}
\begin{proof}[Proof of Theorem \ref{result_wavemap}]
The proof is as above but we additionally assume that $\psi\equiv0$. In this case, we just need the assumptions in Proposition \ref{prop-evol-energy-map} and not the slightly stronger ones in Proposition \ref{prop-evol-energy-spinor}. As the spinor is not involved, we obviously can also remove the spin condition.
\end{proof}

\par\medskip
\textbf{Acknowledgements:}
The first author gratefully acknowledges the support of the Austrian Science Fund (FWF) 
through the START-Project Y963-N35 of Michael Eichmair and
the project P30749-N35 ``Geometric variational problems from string theory''.

The second author wants to thank the Fields Institute for Research in Mathematical Sciences 
for its hospitality during the program on Geometric Analysis.

Both authors thank the Erwin Schrödinger International Institute for Mathematics and
Physics and the organizers of its “Geometry and Relativity” program, where part of this
paper was written.

\bibliographystyle{plain}
\bibliography{mybib}
\end{document}